\documentclass[oneside, a4paper]{amsart}

	\usepackage[a4paper, margin=1in]{geometry}
	\usepackage{amsmath, amsthm, amssymb, latexsym, mathtools, mathrsfs, eucal}
	\usepackage[all,cmtip]{xy}
	\usepackage{dsfont, soul, stackrel, tikz-cd, graphicx, etoolbox, multirow}
	\usepackage{stmaryrd}
	\DeclareGraphicsExtensions{.pdf,.png,.jpg}
	\usepackage{fancyhdr}
	\usepackage[shortlabels]{enumitem}
	\usepackage[pdfmenubar=true, pdfborder={0 0 0 [3 3]}, colorlinks, citecolor=magenta, urlcolor=blue]{hyperref}
    \usepackage{comment}
	\numberwithin{equation}{section}
	\setenumerate{listparindent=\parindent}

	\newtheoremstyle{Mytheorem}%
	{1em}{1em}%
	{\slshape}{}%
	{\bfseries}{.}%
	{ }{}

	\newtheoremstyle{Mydefinition}%
	{1em}{1em}%
	{}{}%
	{\bfseries}{.}%
	{ }{}

	\theoremstyle{Mydefinition}
	\newtheorem{statement}{Statement}[section]
	\newtheorem{definition}[statement]{Definition}
	
	\newtheorem{remark}[statement]{Remark}
	
	\newtheorem{example}[statement]{Example}

	\newtheorem*{comment*}{Comment}

	\theoremstyle{Mytheorem}
	\newtheorem{theorem}[statement]{Theorem}
	\newtheorem{corollary}[statement]{Corollary}
	
	\newtheorem{proposition}[statement]{Proposition}
	\newtheorem{lemma}[statement]{Lemma}

	\newcommand{\G}{GL_2^{+}(\mathbb{Q})}

	\newcommand{\nc}{\newcommand}
	\newcommand{\be}{\begin{eqnarray*}}
	\newcommand{\ee}{\end{eqnarray*}}
	\newcommand{\bea}{\begin{eqnarray}}
	\newcommand{\eea}{\end{eqnarray}}
	\newcommand{\bs}{\begin{split}}
	\newcommand{\es}{\end{split}}
	\newcommand{\bal}{\begin{align}}
	\newcommand{\eal}{\end{align}}

	\nc{\bei}{\begin{itemize}}
	\nc{\eei}{\end{itemize}}
	\nc{\bee}{\begin{enumerate}}
	\nc{\eee}{\end{enumerate}}
	\nc{\bet}{\begin{thm}}
	\nc{\eet}{\end{thm}}
	\nc{\bed}{\begin{defn}}
	\nc{\eed}{\end{defn}}
	\nc{\bel}{\begin{lem}}
	\nc{\eel}{\end{lem}}
	\nc{\bep}{\begin{prop}}
	\nc{\eep}{\end{prop}}
	\nc{\bec}{\begin{corollary}}
	\nc{\eec}{\end{corollary}}
	\nc{\ber}{\begin{rem}}
	\nc{\eer}{\end{rem}}
	\nc{\beex}{\begin{example}}
	\nc{\eeex}{\end{example}}
	\nc{\bpm}{\begin{pmatrix}}
	\nc{\epm}{\end{pmatrix}}
	\nc{\bspm}{\left(\begin{smallmatrix}}
	\nc{\espm}{\end{smallmatrix}\right)}

	\newcommand{\cA}{\mathcal{A}}

	\newcommand{\cF}{\mathcal{F}}
	
	\newcommand{\cH}{\mathcal{H}}

	\newcommand{\cL}{\mathcal{L}}

	\newcommand{\cO}{\mathcal{O}}

	\newcommand{\cT}{\mathcal{T}}
	\newcommand{\cU}{\mathcal{U}}

	\newcommand{\bC}{\mathbb{C}}

	\newcommand{\bM}{\mathbb{M}}

	\newcommand{\bP}{\mathbf{P}}
	
	\newcommand{\bR}{\mathbb{R}}

	\newcommand{\BP}{\mathbf{P}}

	\nc{\frf}{\mathfrak{f}}

	\nc{\frs}{\mathfrak{s}}  
	\nc{\frt}{\mathfrak{t}} 
	\nc{\fru}{\mathfrak{u}}
	\nc{\lsl}{\mathfrak{sl}}
	\nc{\lgl}{\mathfrak{gl}}
	\nc{\upsi}{\underline{\psi}}
	\nc{\uchi}{\underline{\chi}}
	 
	\DeclareMathOperator{\Tr}{Tr}
	\DeclareMathOperator{\PV}{PV}

	\DeclareMathOperator{\Spec}{Spec}

	\newcommand{\lra}{\longrightarrow}    
	
	\nc{\surjto}{\twoheadrightarrow}
	\nc{\ts}{\times}
	\nc{\ds}{\displaystyle}
	\nc{\nd}{\noindent}  
	\nc{\ud}{\underline}
	\nc{\ov}{\overline}
	\nc{\maplra}[1]{\buildrel #1 \over \lra}
	\nc{\mapto}[1]{\buildrel #1 \over \to}
	\nc{\setb}[1]{\{  #1\}}
	\nc{\cHom}{\mathcal{H}om}

\usepackage{MnSymbol}

	\def\a{\alpha}
	\def\b{\beta}

	\def\g{\gamma} \def\G{\Gamma}
	\def\k{\kappa}
	\def\l{\lambda}

	\def\o{\omega} \def\O{\Omega}

	\def\t{\tau}

	\def\U{\Upsilon}

	\def\C{\mathbb{C}}
	
	\def\Z{\mathbb{Z}}

	\def\Ker{\hbox{Ker}\;}
	\def\Im{\hbox{Im}\;}

\newcommand{\Crit}{\operatorname{Crit}}

\title[$L^2$-Hodge theory of Landau-Ginzburg models of smooth projective hypersurfaces]{$L^2$-Hodge theoretic construction of Frobenius manifolds for Calabi-Yau smooth projective hypersurfaces}

\author{Jeehoon Park}
\address{Jeehoon Park: QSMS, Seoul National University, 1 Gwanak-ro, Gwanak-gu, Seoul 08826, South Korea}
\email{jpark.math@gmail.com}

\author{Jaewon Yoo}
\address{Department of Mathematics, POSTECH (Pohang University of Science and Technology), San 31, Hyoja-Dong, Nam-Gu, Pohang, Gyeongbuk, 790-784, South Korea. }
\email{yooj1215@postech.ac.kr}


\begin{document}
\maketitle

\begin{abstract}
We provide a new $L^2$-Hodge theoretic construction of a Frobenius manifold structure on the cohomology of a Calabi-Yau smooth projective hypersurface $V$, using Li-Wen's $L^2$-Hodge theory \cite{LW} of a Landau-Ginzburg model with compact critical locus $V$. 
We also give a precise comparison result between the current construction and Barannikov-Kontsevich's construction \cite{BK} of the Frobenius manifold structure on the cohomology of $V$.

\end{abstract}

\tableofcontents

\section{Introduction}

In \cite{LW}, Si Li and Hao Wen studied the Landau-Ginzburg B-model associated to the triple $(X,\omega_X,W)$ where
\begin{enumerate}
\item
$X$ is a non-compact K\"ahler manifold with complete K\"ahler metric $g$;
\item
$\omega_X$ is a nowhere vanishing holomorphic volume form, i.e. Calabi-Yau form on $X$; 
\item
$W:X\to \bC$ is a holomorphic function with compact critical locus. 
\end{enumerate}

The main result of \cite{LW} is an $L^2$-Hodge theoretic construction of a Frobenius manifold structure on the twisted de Rham cohomology associated with $(X, W)$ under the technical assumptions that $(X, g, \omega_X)$ is a bounded Calabi-Yau geometry (\cite[Definition 3.6]{LW}) and $W$ is strongly elliptic (\cite[Definition 2.6]{LW}).
 In \cite[Section 2.4]{LW}, Li-Wen gave three classes of examples $(X,\omega_X,W)$ satisfying the above technical assumptions: quasi-homogeneous polynomials with an isolated critical point \cite[Subsection 2.4.1]{LW}, convenient Laurent polynomials on $(\C^*)^n$ \cite[Subsection 2.4.3]{LW}, and crepant resolutions of Landau-Ginzburg orbifolds \cite[Subsection 2.4.2]{LW}.
 The construction of Frobenius manifolds for the first and second classes of examples was known before Li-Wen's work \cite{LW}. The construction of Frobenius manifolds in the third class appeared first in \cite{LW}. The goal of the present article is to analyze a specific example in the third class, namely, to compute the twisted de Rham cohomology of a crepant resolution of a Landau-Ginzburg orbifold when the orbifold group is the cyclic group.


Let $n \geq 2$ be any integer. Let $\mu_n$ be the group of $n$-th roots of unity in $\C^*$ and consider the following group action of $\mu_n$ on $\C^n$:
\be
\zeta \cdot (x_1, \cdots, x_n) = (\zeta x_1, \cdots, \zeta x_n), \quad \zeta \in \mu_n, \quad (x_1, \cdots, x_n) \in \C^n.
\ee
Let $f$ be a $\mu_n$ invariant polynomial on $\C^n$ with an isolated singularity at the origin.
We further assume that $f$ is a homogeneous polynomial, which implies that $\deg f =n$ and the zero locus of $f$ defines a Calabi-Yau smooth projective (ample) hypersurface inside the projective space $\BP^{n-1}$:
\be
V_f \subset \BP^{n-1},
\ee
and we use the notation $V:=V_f(\C) \subset \C\BP^{n-1}$ to denote the associated compact K\"ahler Calabi-Yau manifold.

Let $\varpi: X \to \C^n/\mu_n$ be a crepant resolution of $\C^n/\mu_n$. Let $\omega$ be a trivial Calabi-Yau $n$-form $d{x_1} \wedge \cdots \wedge d{x_n}$ on $\C^n$ and denote by $\omega_X$ the pull-back of $\omega$ along $\pi$.
Thus, $(X, \omega_X)$ is a non-compact Calabi-Yau manifold.
Since $\mu_n$ acts freely on $\C^n$ away from the origin, one can apply D. Joyce's result \cite{Joyce} to show that there exists a Ricci-flat ALE (Asymptotically Locally Euclidean) metric $g$ on $X$.
In \cite{LW}, Li-Wen observed the following:
\begin{enumerate}
\item The data $(X, g, \omega_X)$ define a bounded Calabi-Yau geometry in the sense of \cite[Definition 3.6]{LW}.
\item The function $W:=\varpi^*(f)$ is a holomorphic function that satisfies the strongly elliptic condition \cite[Definition 2.6]{LW}. 
\end{enumerate}

\vspace{0.5em}
An explicit description of a crepant resolution $\varpi: X \to \C^n/\mu_n$ and the potential $W=\varpi^*(f)$ is well-known (for example, see \cite[Example 3.5]{Joyce}):
\be
X_{\mathrm{CY}}:=\frac{\C^n \setminus \{0\} \times \C}{\C^*}, \quad \zeta \cdot (x_1, \cdots, x_n, p) = (\zeta x_1, \cdots, \zeta x_n, \zeta^{-n} p),  \zeta \in \C^*, (x_1, \cdots, x_n, p) \in \C^n \setminus \{0\} \times \C,
\ee
and 
\be
\varpi: X_{\mathrm{CY}} \to \C^n/\mu_n, \quad [(x_1, \cdots, x_n, p)] \mapsto [(\xi x_1, \cdots, \xi x_n)]
\ee
where $\xi \in \C$ (unique up to $n$-th root of unity) is given by $\xi^n =p$ (note that $\varpi$ is clearly well-defined),
gives a crepant resolution of $\C^n/\mu_n$; in fact, $X_{\mathrm{CY}}$ is isomorphic to the blow up of $\C^n/\mu_n$ along $\{ 0\}$ such that $\varpi^{-1}({\{0\}})\simeq \BP^{n-1}$.
Moreover, it is easy to check that the holomorphic function
\be
W=W(\ud x, p):= p \cdot f(\ud x): X_{\mathrm{CY}} \to \C
\ee
satisfies $\varpi^* (f) = W(\ud x, p)$ and the critical locus of $W$ is isomorphic to $V=V_f(\C) \subset \C\BP^{n-1}$.
We fix a non-vanishing holomorphic volume form $\O_{\mathrm{CY}}$ on $X_{\mathrm{CY}}$ and a K\"ahler metric as in \cite[Definition 2.16]{LW}.

Let $\cT_{X_{\mathrm{CY}}}$ (respectively, $\cT^*_{X_{\mathrm{CY}}}$) be the holomorphic tangent (respectively, cotangent) bundle on ${X_{\mathrm{CY}}}$. 
Let $\PV^{i,j}({X_{\mathrm{CY}}}):=\cA^{0,j}({X_{\mathrm{CY}}}, \wedge^i \cT_{X_{\mathrm{CY}}})$ denote the space of smooth $(0,j)$-forms valued in $\wedge^i \cT_{X_{\mathrm{CY}}}$.
Consider the isomorphism $\cA^{0,j}({X_{\mathrm{CY}}}, \wedge^i \cT_{X_{\mathrm{CY}}}) \cong \cA^{0,j}({X_{\mathrm{CY}}}, \wedge^{n-i}\cT^*_{X_{\mathrm{CY}}}) $ given by contracting with $\O_{\mathrm{CY}}$, which we denote by $\U_{LW}$:
\be
 \U_{LW}: \bigoplus_{i,j=0}^n \cA^{0,j}({X_{\mathrm{CY}}}, \wedge^i \cT_{X_{\mathrm{CY}}}) \xrightarrow{\cong} \bigoplus_{i,j=0}^n\cA^{0,j}({X_{\mathrm{CY}}}, \wedge^{n-i}\cT^*_{X_{\mathrm{CY}}}).
\ee
This induces an isomorphism $\U_{LW}:H^j({X_{\mathrm{CY}}} , \wedge^i \cT_{X_{\mathrm{CY}}} )\xrightarrow{\simeq} H^{n-i,j}({X_{\mathrm{CY}}})$ for each $0\leq i,j\leq n$ after taking the $\overline{\partial }$-cohomology. 
The de Rham differential decomposes as $d = \partial  +\bar \partial$, where
\be
\partial: \cA^{i,j}({X_{\mathrm{CY}}}) \to \cA^{i+1, j}({X_{\mathrm{CY}}}), \quad \bar\partial: \cA^{i,j}({X_{\mathrm{CY}}}) \to \cA^{i, j+1}({X_{\mathrm{CY}}}).
\ee
We define $\bar \partial: \PV^{i,j} ({X_{\mathrm{CY}}}) \to \PV^{i, j+1}({X_{\mathrm{CY}}})$ and $\partial_{\O_{\mathrm{CY}}}: \PV^{i,j} ({X_{\mathrm{CY}}}) \to \PV^{i-1, j}({X_{\mathrm{CY}}})$ by
\be
\bar\partial(\a) \vdash \O_{\mathrm{CY}} = \bar\partial (\a \vdash \O_{\mathrm{CY}}), \quad \partial_{\O_{\mathrm{CY}}} (\a) \vdash \O_{\mathrm{CY}} = \partial (\a \vdash \O_{\mathrm{CY}}).
\ee
Let
\be \label{PVgrading}
\PV^r({X_{\mathrm{CY}}}):= \bigoplus_{r=j-i}\PV^{i,j}({X_{\mathrm{CY}}}), \quad \cA^r({X_{\mathrm{CY}}}):=\bigoplus_{r=i+j}\cA^{i,j}({X_{\mathrm{CY}}})
\ee
and write $|\a|= a$ whenever $\a \in \PV^a({X_{\mathrm{CY}}})$ (and $|\b|=b$ whenever $\b \in \cA^b({X_{\mathrm{CY}}})$).
Note that $\PV({X_{\mathrm{CY}}})=\oplus_r \PV^r({X_{\mathrm{CY}}})$ has a product structure that comes from the product of holomorphic polyvector fields valued in anti-holomorphic differential forms.
We define a bracket on $\PV({X_{\mathrm{CY}}})$:
\be
\{ \a, \b\} := \partial_{\O_{\mathrm{CY}}} (\a \wedge \b) - \partial_{\O_{\mathrm{CY}}}  (\a) \wedge \b - (-1)^{|\a|} \a \wedge \partial_{\O_{\mathrm{CY}}}  (\b), \quad \a,\b \in \PV({X_{\mathrm{CY}}}).
\ee
Then 
\begin{eqnarray*}
    \U_{LW}(\{W, \a\}) = dW \wedge \U_{LW}(\a),\quad \a \in \PV({X_{\mathrm{CY}}}),
\end{eqnarray*}
and 
\bea \label{lwdgbv}
(\PV({X_{\mathrm{CY}}}), \bar\partial_W:= \bar \partial + \{ W, - \}, \partial_{\O_{\mathrm{CY}}})
\eea
is a dGBV algebra (\cite[Definition 3.2]{LW}); in particular,
$(\PV({X_{\mathrm{CY}}}), \bar \partial_W )$ is a commutative differential graded algebra.
Li-Wen proved \cite[Theorem 3.26]{LW} that there exists a $L^2$-Hodge theoretic construction of a Frobenius manifold structure on the cohomology $H(\mathrm{\PV}({X_{\mathrm{CY}}}), \bar \partial_W)$.

 Our first main result is to compute the cohomology $H(\mathrm{\PV}(X_{\mathrm{CY}}), \bar \partial_W)$ in terms of the cohomology of $V$.

\begin{theorem} \label{mthm}
The cohomology of $(\PV(X_{\mathrm{CY}}), \bar \partial_{W})$ is isomorphic to the cohomology $H(V; \C)$ of the smooth projective hypersurface $V \subset \C\BP^{n-1}$.
\end{theorem}

Using Li-Wen's theorem \cite[Theorem 3.26]{LW}, we obtain the following corollary.
\begin{corollary}\label{mcor}
The $L^2$-Hodge theoretic method of Li-Wen provides a Frobenius manifold structure on $H(\PV(X_{\mathrm{CY}}), \bar \partial_W)\simeq H(V; \C)$.
\end{corollary}

Li-Wen's construction of Frobenius manifolds \cite{LW} uses the $L^2$-Hodge theory of $W$-twisted Sobolev spaces.
On the other hand, there is another Frobenius manifold structure on $H(V;\C)$ due to Barannikov-Kontsevich \cite{BK}. Similar to \eqref{lwdgbv}, one has a dGBV algebra 
\bea \label{bkdgbv}
(\PV(V), \bar\partial, \partial_{\O_V})
\eea
on the Calabi-Yau compact manifold $V$, where $\O_V$ is a non-vanishing holomorphic volume form; this is due to Barannikov-Kontsevich \cite{BK}; they proved that  
there exists a formal Frobenius manifold structure on the cohomology $H(\mathrm{\PV}(V), \bar \partial)\simeq H(V;\C)$. 
Theorem \ref{mthm} leads to a natural question:
\be
\text{How do we compare these two Frobenius manifold structures on $H(V;\C)$?}
\ee

Both constructions of Frobenius manifolds \cite{BK} and \cite{LW} are based on the same method, that is, the dGBV algebra with \textit{an integral} satisfying suitable conditions provides a Frobenius manifold structure on the cohomology of the dGBV algebra; see \cite[Section 6]{Ma03} for details.
We have two dGBV algebras, namely, the LW (Li-Wen) dGBV algebra (cf. \eqref{lwdgbv} and see \cite[Definition 3.4]{LW} for the definition of $\PV_{W,\infty}(X_{\mathrm{CY}})$)
\begin{eqnarray}\label{LWD}
    (\cA_{LW},\bar\partial_W, \partial_{\mathrm{CY}}) := (\PV_{W,\infty}(X_{\mathrm{CY}}), \bar\partial_W, \partial_{\O_{\mathrm{CY}}})
\end{eqnarray}
and the BK (Barannikov-Kontsevich) dGBV algebra (cf. \eqref{bkdgbv})
\begin{eqnarray}\label{BKD}
    (\cA_{BK}, \bar\partial, \partial_{\O_f})
     :=(\PV(V),\bar\partial, \partial_{\O_f})
\end{eqnarray}
where $\O_f$ is given in \eqref{volume}.
These dGBV algebras give us two Frobenius algebras, namely, the LW Frobenius algebra 
$$
(\O_{LW}, [1], K_{LW}^{(0)})=(H(\cA_{LW},\bar\partial_W),[1], K_{LW}^{(0)} )
$$
and the BK Frobenius algebra 
$$
(\O_{BK},[1],K_{BK}^{(0)})=(H(\cA_{BK},\bar\partial), [1], K_{BK}^{(0)} );
$$ 
see Definitions \ref{LWF} and \ref{BKF} for details.
Since two dGBV structures are unlikely to be isomorphic, it seems to be difficult to compare two Frobenius manifold structures directly; see Subsection \ref{sec3.3} for more comments. But we show that their Frobenius algebras are isomorphic up to a nonzero constant:

\begin{theorem}\label{theoremtwo}
There is a ring isomorphism
    \begin{eqnarray*}
        \Phi:\O_{LW} \to \O_{BK}
    \end{eqnarray*}
    such that there exists a non-zero constant $c_\Phi \in \C$ satisfying
    \begin{eqnarray*}
        c_\Phi K_{LW}^{(0)}(\a, \b) = K_{BK}^{(0)}(\Phi(\a), \Phi(\b)), \quad \a, \b \in \O_{LW}.
    \end{eqnarray*}
\end{theorem}
{
A key new input for its proof is an explicit computation of the trace $\Tr_{LW}$(Lemma \ref{integrationcompute}) and the residue pairing $K_{LW}^{(0)}$ (Proposition \ref{proph} and Proposition \ref{ringp}) using the homotopy operator $T_\rho:\PV(X_{\mathrm{CY}}) \to \PV_{W,\infty}(X_{\mathrm{CY}})$ in \cite{LLS}, which serves as a cohomological inverse operator to the embedding quasi-isomorphism $\iota_2: \PV_{W,\infty}(X_{\mathrm{CY}}) \to \PV(X_{\mathrm{CY}})$ (see \cite[Theorem 2.34]{LW}).

 The ring isomorphism $\Phi$ is not induced from a cochain map $(\cA_{LW}, \bar \partial_W) \to (\cA_{BK}, \bar \partial)$; it seems unlikely that a cochain map $(\cA_{LW}, \bar \partial_W) \to (\cA_{BK}, \bar \partial)$, which induces a vector space isomorphism (not to mention a ring isomorphism) between $\O_{LW}$ and $\O_{BK}$, would exist.

\vspace{0.5em}

\textbf{Structure of the paper:} 
In Section \ref{sec2}, we compute the cohomology of the commutative differential graded algebra $(\PV(X_{\mathrm{CY}}), \bar\partial_W)$ and show that it is isomorphic to the cohomology of the projective smooth hyperspace $V=V_f(\C)$: we prove Theorem \ref{mthm}.
To this end, we introduce two spectral sequences $E^{p,q}$ and $\tilde E^{p,q}$, namely, analytic and algebraic spectral sequences. In Subsection \ref{sec2.1}, we compute the $E_2$-page of the algebraic spectral sequence $\tilde E^{p,q}$. In Subsection \ref{sec2.2}, we prove that the $E_2$-pages of two spectral sequences are the same. In Subsection \ref{sec2.3}, we compute the higher terms in the spectral sequences to conclude Theorem \ref{mthm}.

In Section \ref{sec3}, we provide a comparison result between Li-Wen's Frobenius algebra and Barannikov-Kontsevich's Frobenius algebra: we prove Theorem \ref{theoremtwo}. In Subsection \ref{sec3.1}, we introduce the LW Frobenius algebra and the BK Frobenius algebra. In Subsection \ref{sec3.2} we show that they are isomorphic up to constant factor on the inner product. We give some remark on the comparison of Frobenius manifold structures in Subsection \ref{sec3.3}. 
\vspace{0.5em}

\textbf{Acknowledgement}:
Jeehoon Park thanks Philsang Yoo for a useful and valuable discussion on the subject. Jeehoon Park was supported by the National Research Foundation of Korea (NRF-2021R1A2C1006696) and the National Research Foundation of Korea (NRF) grant funded by the Korea government (MSIT) (No.2020R1A5A1016126).

\section{Computation of the cohomology}\label{sec2}

We define the Jacobian ring $R(W)$ as
\be
R(W):=\frac{\bC[\ud{x},p]}{J(W)}=\frac{\bC[x_1, \cdots, x_n,p]}{J(W)}
\ee
where $J(W)=(\frac{\partial W}{\partial x_1},\ldots,\frac{\partial W}{\partial x_n}, \frac{\partial W}{\partial p})$ is the Jacobian ideal of $\bC[\ud x,p]$.
We define an additive grading, which we call charge $ch$:
\begin{eqnarray*}
    \mathrm{ch}(x_1)=1, \ldots, \mathrm{ch}(x_n)=1, \quad \mathrm{ch}(p)=-n.
\end{eqnarray*}
This extends to the $\C[\ud x, p]$-module $\O_{\C[\ud x, p]/\C}^1$ of K\"ahler differentials:
\begin{align*}
\mathrm{ch}(dx_1)=\cdots=\mathrm{ch}(dx_n)=1, \quad
\mathrm{ch}(dp)=-n.
\end{align*}
For each integer $m$, let
\be
R(W)_m := \frac{\C[\ud x,p]_{m}}{\C[\ud x,p]_{m} \cap J(W)}
\ee
where $\C[\ud x,p]_{m}$ is the space of homogeneous elements of charge $m$ in $\C[\ud x,p]$.
We also define the Jacobian ring $R(f)$ as
\be
R(f):=\frac{\bC[\ud{x}]}{J(f)}=\frac{\bC[x_1, \cdots, x_n]}{J(f)}
\ee
where $J(f)=(\frac{\partial f}{\partial x_1},\ldots,\frac{\partial f}{\partial x_n})$ is the Jacobian ideal of $\bC[\ud x]$.
For each integer $m$, let
\be
R(f)_m := \frac{\C[\ud x]_m}{\C[\ud x]_m \cap J(f)}
\ee
where $\C[\ud x]_m$ is the homogeneous degree $m$ component of $\C[\ud x]$. Note that there is a ring isomorphism 
\be
\bigoplus_{k=0}^\infty R(f)_{kn} \xrightarrow{\sim} R(W)_0, \quad [u(\ud x)] \mapsto [p^{k}u(\ud x)]
\ee
where $u(\ud x) \in \C[\ud x]_{kn}$ is a homogeneous polynomial of degree $kn$.

The weak Lefschetz theorem and P. Griffith's theorem \cite{Gr69} says the following:
\bea \label{cohv}
H^{r}(V_f;\bC) \cong \begin{cases}
\bC, & \textrm{if $0\leq {r}\leq 2(n-2)$, $q\neq n-2$, ${r}$ is even,} \\
R(W)_{0}, & \textrm{if ${r}=n-2$, $n-2$ is odd,}\\
R(W)_{0} \oplus \bC, & \textrm{if ${r}=n-2$, $n-2$ is even,}\\
0, & \textrm{otherwise.}
\end{cases}
\eea
In particular, $\dim_\C H(V_f;\C) = \dim_\C R(W)_0 + n-1$.
In the previous section, we considered $X_{\mathrm{CY}}$ as a complex manifold. Here we need to view $X_{\mathrm{CY}}$ as an algebraic variety over $\bC$, endowed with the Zariski topology and the structure sheaf $\cO_{X_{\mathrm{CY}}}$ whose value on an open subset $U\subset X_{\mathrm{CY}}$ is an algebraic (in other words, polynomial) function on $U$. The complex manifold we have denoted by $X_{\mathrm{CY}}$ will be denoted by $X_{\mathrm{CY}}^{an}$, which is the analytification of $X_{\mathrm{CY}}$; see \cite[p. 7]{SerreGAGA}. 
Then $X_{\mathrm{CY}}^{an}$ is also a ringed space whose topology is the Euclidean topology, and the stalk $(\cO_{X_{\mathrm{CY}}}^{an})_{[\ud{x},p]}$ at $[\ud{x},p]\in X_{\mathrm{CY}}^{an}$ consists of germs of holomorphic functions at $[\ud{x},p]$. Consequently, we will distinguish the projective variety $\bP^{n-1}_{\bC}$ and the compact complex manifold $(\bP^{n-1}_{\bC})^{an}=\bC\bP^{n-1}$. 
We are given a ringed space morphism $h:(X_{\mathrm{CY}}^{an},\cO_{X_{\mathrm{CY}}}^{an})\to (X_{\mathrm{CY}},\cO_{X_{\mathrm{CY}}})$, which maps a point in $X_{\mathrm{CY}}^{an}$ to its corresponding maximal ideal in $X_{\mathrm{CY}}$. There is an analytification functor $(-)^{an}:QCoh(\cO_{X_{\mathrm{CY}}})\to QCoh(\cO_{X_{\mathrm{CY}}}^{an})$ between the category of quasi-coherent sheaves on $(X_{\mathrm{CY}},\cO_{X_{\mathrm{CY}}})$ and $(X_{\mathrm{CY}}^{an},\cO_{X_{\mathrm{CY}}}^{an})$ defined as
$$
\cF^{an} = h^{-1} \cF \otimes_{h^{-1}\cO_{X_{\mathrm{CY}}}}\cO_{X_{\mathrm{CY}}}^{an}.
$$
See \cite[p. 16]{SerreGAGA} for more details.

Using the spectral sequence of the double complex where the vertical differential is $\overline{\partial}$ and the horizontal differential is $dW\wedge$, we have a spectral sequence whose 0th, 1st, and 2nd pages are
\begin{align*}
E_0^{{r},s} &= \cA^{{r},s}(X_{\mathrm{CY}}^{an}),\\
E_1^{{r},s} &= H^s(\cA^{{r},\bullet}(X_{\mathrm{CY}}^{an}),\overline{\partial})\cong H^s(X_{\mathrm{CY}}^{an},(\Omega_{X_{\mathrm{CY}}}^{{r}})^{an}),\\
E_2^{{r},s} &= H^{r}(H^s(\cA^{\bullet,\bullet}(X_{\mathrm{CY}}^{an}),\overline{\partial}),dW)\cong H^{r}(H^s(X_{\mathrm{CY}}^{an},(\Omega_{X_{\mathrm{CY}}}^{\bullet})^{an}),dW),
\end{align*}
where $\O_{X_{\mathrm{CY}}}$ is the sheaf of algebraic differential forms on $X_{\mathrm{CY}}$ and $(\O_{X_{\mathrm{CY}}})^{an}$ is the sheaf of holomorphic differential forms on $X_{\mathrm{CY}}^{an}$.
\begin{comment}
 By Dolbeault's theorem, we have
$$
H^q(\cA^{p,\bullet}(X_{\mathrm{CY}}^{an}),\overline{\partial}) \cong H^q(X_{\mathrm{CY}}^{an},(\Omega_{X_{\mathrm{CY}}}^{p})^{an}).
$$
So we have
\begin{align*}
E_1^{p,q} &\cong H^q(X_{\mathrm{CY}}^{an},(\Omega_{X_{\mathrm{CY}}}^{p})^{an}),\\
E_2^{p,q} &\cong H^p(H^q(X_{\mathrm{CY}}^{an},(\Omega_{X_{\mathrm{CY}}}^{\bullet})^{an}),dW).
\end{align*}
\end{comment}
The cohomology is easier to manipulate in algebraic setting; we have another spectral sequence of the double complex where the vertical differential is the \v{C}ech differential, and the horizontal differential is defined as $dW\wedge$. We define $\cU:=\{U_1,\ldots,U_n\}$  as an open cover of $X_{\mathrm{CY}}$, where each $U_j$ is defined as the preimage of the affine open subset defined by $x_j \neq 0$ inside $\bP^{n-1}_{\bC}$ along the natural projection map $\pi:X_{\mathrm{CY}}\to \bP_{\bC}^{n-1}$. Their intersections $U_{i_0,\ldots,i_{r}}=U_{i_0}\cap \cdots\cap U_{i_{r}}$ are all Stein. Then the 0th page of the spectral sequence is defined by the \v{C}ech complex:
\begin{align*}
\widetilde{E}_0^{{r},s} &= \check{C}^{s}(\cU,\Omega^{r}_{X_{\mathrm{CY}}}),\\
\widetilde{E}_1^{{r},s} &= H^s(X_{\mathrm{CY}},\Omega^{r}_{X_{\mathrm{CY}}}),\\
\widetilde{E}_2^{{r},s} &= H^q(H^s(X_{\mathrm{CY}},\Omega^{\bullet}_{X_{\mathrm{CY}}}),dW).
\end{align*}
 We compute the second page of this spectral sequence in Subsection \ref{sec2.1}. 
\begin{lemma}
$$
\widetilde{E}_2^{{r},s}=H^{r}(H^s(X_{\mathrm{CY}},\Omega_{X_{\mathrm{CY}}}^{\bullet}),dW)\simeq\begin{cases}
\bC, & \textrm{if $1\leq {r}=s\leq n-1$,}\\
R(W)_{0}, & \textrm{if $s=0, {r}=n$}, \\
0, &\textrm{otherwise.}
\end{cases}
$$
\label{finitedimension}
\end{lemma}

The cohomologies of the sheaf of holomorphic differential forms $(\Omega^{s}_{X_{\mathrm{CY}}})^{an}$ and the sheaf of algebraic differential forms $\Omega^{s}_{X_{\mathrm{CY}}}$ are not the same in general, because $X_{\mathrm{CY}}$ is not compact. We use Lemma \ref{finitedimension} (the finite-dimensionality of $\widetilde{E}_2^{{r},s}$) to prove the following lemma, which connects algebraic cohomology with analytic cohomology, i.e. $E_2^{{r},s}\cong \widetilde{E}_2^{{r},s}$:
\begin{lemma}
\label{lemma5.2}
For each ${r},s$, we have an isomorphism
$$H^{r}(H^s(X_{\mathrm{CY}}^{an},(\Omega_{X_{\mathrm{CY}}}^{\bullet})^{an}),dW) \cong H^{r}(H^s(X_{\mathrm{CY}},\Omega_{X_{\mathrm{CY}}}^{\bullet}),dW).$$
\end{lemma}

In Subsection \ref{sec2.3}, we show that all higher differentials $d_k^{{r},s}:E_k^{{r},s} \to E_k^{{r}+k,s-k+1}$ with $k>2$ are trivial; see Lemma \ref{higherdiffan}. Hence $E_2^{{r},s}\simeq E_{\infty}^{{r},s}$
so that we obtain the following theorem, which proves Theorem \ref{mthm} by using \eqref{cohv}.
\begin{theorem}
\label{cohomologycmp}
$$
H^{{r}}(PV^\bullet(X_{\mathrm{CY}}),\overline{\partial}_W) \cong \begin{cases}
\bC, & \textrm{if ${r}=-n+2+2k\neq 0$, $0\leq k \leq n-2$}\\
R(W)_{0}, & \textrm{if ${r}=0$ and $n-2$ is odd,}\\
R(W)_{0}\oplus \bC, & \textrm{if ${r}=0$ and $n-2$ is even,}\\
0, & \textrm{otherwise.}
\end{cases}
$$
\end{theorem}

\subsection{Computation of the second page of the algebraic spectral sequence}
\label{sec2.1}

We prove Lemma  \ref{finitedimension}. Using the Leray spectral sequence for the projection map $\pi:X_{\mathrm{CY}} \to \BP^{n-1}_\C$, we have the following lemma.
\begin{proposition}
\label{leray}
For each ${r},s \geq 0$, we have an isomorphism
\begin{align*}
H^s(X_{\mathrm{CY}},\Omega^{{r}}_{X_{\mathrm{CY}}}) \simeq H^s(\bP^{n-1}_{\bC},\pi_* \Omega^{{r}}_{X_{\mathrm{CY}}}).
\end{align*}
\end{proposition}
\begin{proof}
It follows directly from the fact that $\pi$ is an affine morphism and \cite[Lemma 01F4]{Stacks}.
\end{proof}

\begin{proposition}[Euler's sequence] \label{lone}
We have a short exact sequence of $\cO_{X_{\mathrm{CY}}}$-modules
\begin{align*}
0\to \Omega_{X_{\mathrm{CY}}} \to \{\cO_{X_{\mathrm{CY}}}(-1)\}^{\oplus n}\oplus\cO_{X_{\mathrm{CY}}}(n)\to \cO_{X_{\mathrm{CY}}} \to 0,\\
0\to \Omega_{X_{\mathrm{CY}}}^{{r}} \to \bigwedge^{r} \left(\{\cO_{X_{\mathrm{CY}}}(-1)\}^{\oplus n}\oplus \cO_{X_{\mathrm{CY}}}(n)\right)\to \Omega_{X_{\mathrm{CY}}}^{{r}-1} \to 0.
\end{align*}
\end{proposition}
\begin{proof}
The second short exact sequence is induced by taking the wedge product to the first, together with the fact that $\Omega_{X_{\mathrm{CY}}}$ is a locally free sheaf. On $U_j$ defined by $x_j\neq 0$, define an affine coordinate $X_k=x_k/x_j$ for $k=1,\ldots,n$ and $P = x_j^{n} p$. Then $\Omega_{X_{\mathrm{CY}}}(U_j)$ is a free $\cO_{X_{\mathrm{CY}}}(U_j)$-module of rank $n$ whose bases is $\{dX_1,\ldots,dX_{j-1},dX_{j+1},\ldots,dX_n,dP\}$. We identify
$$
\{\cO_{X_{\mathrm{CY}}}(-1)\}^{\oplus n}\oplus \cO_{X_{\mathrm{CY}}}(n) \simeq \cO_{X_{\mathrm{CY}}}(n) dp \oplus \bigoplus_{i=1}^n \cO_{X_{\mathrm{CY}}}(-1) dx_i.
$$ 
so that
$$
f_1 dx_1 + \cdots + f_n dx_n + g dp 
$$
is an element $(f_1,\cdots,f_n,g)\in \left(\{\cO_{X_{\mathrm{CY}}}(-1)\}^{\oplus n}\oplus \cO_{X_{\mathrm{CY}}}(n)\right)(U_j)$.
The map
$$
\Omega_{X_{\mathrm{CY}}}(U_j) \to \left(\{\cO_{X_{\mathrm{CY}}}(-1)\}^{\oplus n}\oplus\cO_{X_{\mathrm{CY}}}(n)\right)(U_j)
$$
is then defined as
$$
dX_k \mapsto \frac{1}{x_j} dx_k - \frac{x_k}{x_j^2} dx_j, \quad dP \mapsto x_j^{n}dp + n x_j^{n-1}p dx_j,
$$
for $k=1,\ldots,n$. The map
$$
\left(\{\cO_{X_{\mathrm{CY}}}(-1)\}^{\oplus n}\oplus \cO_{X_{\mathrm{CY}}}(n)\right)(U_j)\to \cO_{X_{\mathrm{CY}}}(U_j)
$$
is defined as
\begin{align}
\label{eulersequencesurj}
dx_k \mapsto x_k, \quad dp \mapsto -n p,
\end{align}
for $k=1,\ldots,n$.  
One can check that these maps induce morphisms between sheaves of $\cO_{X_{\mathrm{CY}}}$-modules and define the desired exact sequence.
\end{proof}

Since $\pi$ is an affine morphism, we have the following exact sequence after taking the direct image functor $\pi_*$:
\begin{eqnarray*}
    0\to \pi_*\Omega_{X_{\mathrm{CY}}}^{{r}} \to \bigwedge^{r} \left(\{\pi_*\cO_{X_{\mathrm{CY}}}(-1)\}^{\oplus n}\oplus \pi_*\cO_{X_{\mathrm{CY}}}(n)\right)\to \pi_*\Omega_{X_{\mathrm{CY}}}^{{r}-1} \to 0.
\end{eqnarray*}
The sheaf $\pi_* \cO_{X_{\mathrm{CY}}}(n)$ over $\bP^{n-1}_{\bC}$ can be computed as follows.

\begin{lemma}
For any $d \in \Z$, we have
$$
\pi_* \cO_{X_{\mathrm{CY}}}(d) \cong \bigoplus_{k\geq 0} \cO_{\bP_{\bC}^{n-1}}\left(nk+d \right).
$$
\label{lemma5.6}
\end{lemma}
\begin{proof}
We compute $\pi_* \cO_{X_{\mathrm{CY}}}(d)$ on the affine open subset $V_j$ of $\bP_{\bC}^{n-1}$ defined by $x_j\neq 0$: 
\begin{align*}
\pi_* \cO_{X_{\mathrm{CY}}}(d) (V_j) &=  \cO_{X_{\mathrm{CY}}}(d)(\pi^{-1}(V_j)), \\
&= \bC\left[x_1,\ldots,x_n,p,\frac{1}{x_j}\right]_{d}, \\
&\cong \bigoplus_{k\geq 0} \bC\left[x_1,\ldots,x_n,\frac{1}{x_j}\right]_{nk+d} p^k, \\
&\cong \bigoplus_{k\geq 0} \cO_{\bP_{\bC}^{n-1}}\left(nk+d \right)(V_j).
\end{align*}
\end{proof}

\begin{proposition}\cite[{Lemma 01XT}]{Stacks}\label{lthree}
Let $n$ be a positive integer. We have
$$
H^s (\bP^{n-1}_{\bC}, \cO_{\bP_{\bC}^{n-1}}(d)) = \begin{cases}
\bC[\ud x]_{d} & \textrm{if } s=0, d\geq 0, \\
\left(\frac{1}{x_1\cdots x_n}\bC\left[\frac{1}{x_1},\ldots,\frac{1}{x_n}\right]\right)_{d} & \textrm{if }  s=n-1, d\leq  -n, \\
0 & \textrm{otherwise}. 
\end{cases}
$$
\label{lemma5.7}
\end{proposition}

\begin{proof}[Proof of Lemma \ref{finitedimension}]
By Proposition \ref{lthree}, Lemma \ref{lemma5.6}, and Proposition \ref{leray}, the $s$th sheaf cohomology of the middle term is nonzero only if either $s=0$ or $(s,{r})=(n-1,n)$:
\begin{align*}
H^s\left(X_{\mathrm{CY}},\bigwedge^{r} \left(\{\cO_{X_{\mathrm{CY}}}(-1)\}^{\oplus n}\oplus \cO_{X_{\mathrm{CY}}}(n)\right)\right) \simeq 
\begin{cases}
\left(\Omega_{\bC[\ud{x},p]/\bC}^{r}\right)_{0} & \textrm{if $s=0$,} \\
\bC & \textrm{if $s=n-1$, ${r}=n$,} \\
0 & \textrm{otherwise.}
\end{cases}
\end{align*}

Applying the sheaf cohomology to the Euler sequences (Proposition \ref{lone}) for ${r}=0$, we have
$$
0\to H^0(X_{\mathrm{CY}},\Omega_{X_{\mathrm{CY}}}^{{r}})\to(\Omega_{\bC[\ud{x},p]/\bC}^{r})_{0}\to H^0(X_{\mathrm{CY}},\Omega_{X_{\mathrm{CY}}}^{{r}-1})\to H^1(X_{\mathrm{CY}},\Omega_{X_{\mathrm{CY}}}^{{r}}) \to 0.
$$
When ${r}=n$, we have
\begin{align*}
0\to H^{n-2}(X_{\mathrm{CY}},\Omega_{X_{\mathrm{CY}}}^{n-1})\to H^{n-1}(X_{\mathrm{CY}},\Omega_{X_{\mathrm{CY}}}^{n}) \to \bC \to H^{n-1}(X_{\mathrm{CY}},\Omega_{X_{\mathrm{CY}}}^{n-1}) \to H^{n}(X_{\mathrm{CY}},\Omega_{X_{\mathrm{CY}}}^{n}) \to 0.
\end{align*}
If $s\neq 0$, $(r,s)\neq (n-1,n-1)$, and $(r,s)\neq (n-1,n-2)$, then we have
\begin{align}
H^s(X_{\mathrm{CY}},\Omega_{X_{\mathrm{CY}}}^{{r}}) \cong H^{s+1}(X_{\mathrm{CY}},\Omega_{X_{\mathrm{CY}}}^{{r}+1}).
\label{diagiso}
\end{align}
If ${r}\geq n+1$, then $\Omega_{X_{\mathrm{CY}}}^{{r}}=0$, so ${r}$ large enough, $H^s(X_{\mathrm{CY}},\Omega^{{r}}_{X_{\mathrm{CY}}})=0$. Thus, we have the following computation:
\begin{align}
\tilde E_1^{{r},s}=H^s(X_{\mathrm{CY}},\Omega_{X_{\mathrm{CY}}}^{{r}}) = \begin{cases}
\bC & \textrm{if ${r}=s$, $1\leq {r} \leq n-1,$} \\
H^0(X_{\mathrm{CY}},\Omega_{X_{\mathrm{CY}}}^{{r}}) & \textrm{if $s=0$, $0\leq {r} \leq n$,} \\
0 & \textrm{otherwise}.
\end{cases}
\label{algcohpage1}
\end{align}


In order to take the cohomology along $dW$, we only need to compute the cohomology of the chain complex $(H^0(X_{\mathrm{CY}}, \O_{X_{\mathrm{CY}}}^\bullet, dW))\simeq(H^0(\bP^{n-1},\pi_*\Omega^\bullet_{X_{\mathrm{CY}}}),dW)$. 
By Proposition \ref{isolated}, we have
\begin{align}
H^s(H^0(X_{\mathrm{CY}},\Omega^{\bullet}_{X_{\mathrm{CY}}}),dW) \cong \begin{cases}
\bigoplus_{k=0}^\infty R(f)_{kn}, &\textrm{if $s=n$,} \\
0, & \textrm{otherwise.}
\end{cases}
\label{secondpage}
\end{align}
since $f\in \bC[\ud{x}]$ has isolated hypersurface singularity at $0$ and the corresponding Koszul complex is concentrated in the top degree. This finishes the proof of Lemma \ref{finitedimension}
\end{proof}

\begin{proposition}\label{isolated}
    There is an isomorphism of chain complexes between $(H^0(\bP^{n-1},\pi_*\Omega^\bullet_{X_{\mathrm{CY}}}),dW)$ and the following Koszul subcomplex associated to $df$
$$
0\to \bigoplus_{k=0}^\infty \bC[\ud{x}]_{kn} \stackrel{\wedge df}{\longrightarrow} \bigoplus_{k=0}^\infty\left(\Omega_{\bC[\ud{x}]}^1\right)_{kn} \stackrel{\wedge df}{\longrightarrow}  \cdots \stackrel{\wedge df}{\longrightarrow} \bigoplus_{k=0}^\infty\left(\Omega_{\bC[\ud{x}]}^n\right)_{kn} \to 0.
$$
\end{proposition}
\begin{proof}
We explicitly construct a map $\bigoplus_{k=0}^\infty(\Omega^s_{\bC[\ud{x}]})_{kn} \to H^0(\bP^{n-1},\pi_*\Omega_{X_{\mathrm{CY}}}^s)$. Suppose the element $\omega=fdx_{j_1}\wedge\cdots\wedge dx_{j_s}\in \bigoplus_{k=0}^\infty\left(\Omega_{\bC[\ud{x}]}^s\right)_{kn}$ is given. Let $f_{kn-s}$ be the homogeneous component of degree $kn-s$ of $f$. We let $F(\ud{x},p)\in \bC[\ud{x},p]_{-s}$ be the homogenization of $f$, i.e.
$$
F:= \sum_{k=0}^\infty p^k f_{kn-s}.
$$
Unless $s=0$, we have $f_{-s}=0$ for $k=0$. Thus, $F\in (p)\cap \bC[\ud{x},p]_{-s}$ when $s\neq 0$.

Then $\omega$ is mapped to
$$
\omega'=Fdx_{j_1}\wedge\cdots\wedge dx_{j_s}+\frac{1}{np}\sum_{t=1}^s x_{j_t}F dx_{j_1}\wedge\cdots\wedge dx_{j_{t-1}}\wedge dp\wedge dx_{j_{t+1}}\wedge \cdots \wedge dx_{j_{s}}.
$$ 
$\omega'$ is a globally defined section in $\pi_*\Omega_{X_{\mathrm{CY}}}^s$ over $\bP^{n-1}$, since it is in the kernel of the morphism in Proposition \ref{lone},
$$
\bigwedge^s \left(\{\cO_{X_{\mathrm{CY}}}(-1)\}^{\oplus n}\oplus \cO_{X_{\mathrm{CY}}}(n)\right)\to \Omega_{X_{\mathrm{CY}}}^{s-1} \to 0.
$$
We want to show that $dW\wedge \omega'$ is same as the image of $df\wedge \omega$.
\begin{align*}
df\wedge \omega =& df \wedge\left(fdx_{j_1}\wedge\cdots\wedge dx_{j_s}\right) \\
=&\sum_{\substack{1\leq l \leq n\\l\notin\{j_1,\ldots,j_s\}}} f\frac{\partial f}{\partial x_l} dx_l \wedge dx_{j_1}\wedge\cdots\wedge dx_{j_s}= \sum_{\substack{1\leq l \leq n\\l\notin\{j_1,\ldots,j_s\}}}  \sum_{k=0}^\infty f_{kn-s}\frac{\partial f}{\partial x_l} dx_l \wedge dx_{j_1}\wedge\cdots\wedge dx_{j_s}, \\
\mapsto&\sum_{\substack{1\leq l \leq n\\l\notin\{j_1,\ldots,j_s\}}} \sum_{k=0}^\infty \left( p^{k+1}f_{kn-s}  \frac{\partial f}{\partial x_l} dx_l \wedge dx_{j_1}\wedge\cdots\wedge dx_{j_s}\right. \\
&\qquad \left.+\frac{1}{n} \sum_{t=1}^s x_{j_t}p^k f_{kn-s}\frac{\partial f}{\partial x_l} dx_l\wedge dx_{j_1}\wedge\cdots\wedge dx_{j_{t-1}}\wedge dp\wedge dx_{j_{t+1}}\wedge \cdots \wedge dx_{j_{s}}\right.\\
&\qquad  \left.+\frac{1}{n}x_l p^k f_{kn-s}\frac{\partial f}{\partial x_l} dp\wedge dx_{j_1}\wedge\cdots\wedge dx_{j_s}\right), \\
=& \sum_{\substack{1\leq l \leq n\\l\notin\{j_1,\ldots,j_s\}}} \left( F  \frac{\partial W}{\partial x_l} dx_l \wedge dx_{j_1}\wedge\cdots\wedge dx_{j_s}\right. \\
&\qquad \left.+\frac{1}{np} \sum_{t=1}^s x_{j_t}F\frac{\partial W}{\partial x_l} dx_l\wedge dx_{j_1}\wedge\cdots\wedge dx_{j_{t-1}}\wedge dp\wedge dx_{j_{t+1}}\wedge \cdots \wedge dx_{j_{s}}\right.\\
&\qquad \left.+\frac{1}{np} F x_l\frac{\partial W}{\partial x_l} dp\wedge dx_{j_1}\wedge\cdots\wedge dx_{j_s}\right).
\end{align*}
For each $t=1,\ldots,s$, we have the identity
$$
x_{j_t}F\frac{\partial W}{\partial x_{j_t}} dx_{j_t} \wedge dx_{j_1}\wedge\cdots\wedge dx_{j_{t-1}}\wedge dp\wedge dx_{j_{t+1}}\wedge \cdots \wedge dx_{j_{s}} = - x_{j_t}F\frac{\partial W}{\partial x_{j_t}} dp \wedge dx_{j_1}\wedge\cdots \wedge dx_{j_{s}}.
$$
Thus,
\begin{align*}
&\sum_{\substack{1\leq l \leq n\\l\notin\{j_1,\ldots,j_s\}}} \Bigg(\sum_{t=1}^s x_{j_t}F\frac{\partial W}{\partial x_l} dx_l\wedge dx_{j_1}\wedge\cdots\wedge dx_{j_{t-1}}\wedge dp\wedge dx_{j_{t+1}}\wedge \cdots \wedge dx_{j_{s}}\\
&\qquad + F x_l\frac{\partial W}{\partial x_l} dp\wedge dx_{j_1}\wedge\cdots\wedge dx_{j_s}\Bigg)\\
=&\sum_{\substack{1\leq l \leq n\\l\notin\{j_1,\ldots,j_s\}}} \sum_{t=1}^s x_{j_t}F\frac{\partial W}{\partial x_l} dx_l\wedge dx_{j_1}\wedge\cdots\wedge dx_{j_{t-1}}\wedge dp\wedge dx_{j_{t+1}}\wedge \cdots \wedge dx_{j_{s}}\\
&\qquad +\sum_{\substack{1\leq l \leq n\\l\notin\{j_1,\ldots,j_s\}}} F x_l\frac{\partial W}{\partial x_l} dp\wedge dx_{j_1}\wedge\cdots\wedge dx_{j_s}\\
&\qquad +\sum_{t=1}^s x_{j_t}F\frac{\partial W}{\partial x_{j_t}} dx_{j_t} \wedge dx_{j_1}\wedge\cdots\wedge dx_{j_{t-1}}\wedge dp\wedge dx_{j_{t+1}}\wedge \cdots \wedge dx_{j_{s}} \\
&\qquad + \sum_{t=1}^s F x_{j_t}\frac{\partial W}{\partial x_{j_t}} dp \wedge dx_{j_1}\wedge\cdots \wedge dx_{j_{s}} \\
=& \sum_{t=1}^s x_{j_t}F \left(\sum_{l=1}^n \frac{\partial W}{\partial x_l} dx_l\right)\wedge dx_{j_1}\wedge\cdots\wedge dx_{j_{t-1}}\wedge dp\wedge dx_{j_{t+1}}\wedge \cdots \wedge dx_{j_{s}} \\
&\qquad + F\left(\sum_{l=1}^n x_l \frac{\partial W}{\partial x_l} \right)dp\wedge dx_{j_1}\wedge \cdots \wedge dx_{j_s}.
\end{align*}
By the homogeneity of $W$, we have the equality
$$
\sum_{l=1}^n x_l \frac{\partial W}{\partial x_l} = nW.
$$
Thus, the image of $df\wedge \omega$ is
\begin{align*}
&\sum_{\substack{1\leq l \leq n\\l\notin\{j_1,\ldots,j_s\}}} F  \frac{\partial W}{\partial x_l} dx_l \wedge dx_{j_1}\wedge\cdots\wedge dx_{j_s}\\
&\qquad +\frac{1}{np}\sum_{t=1}^s x_{j_t}F \left(\sum_{l=1}^n \frac{\partial W}{\partial x_l} dx_l\right)\wedge dx_{j_1}\wedge\cdots\wedge dx_{j_{t-1}}\wedge dp\wedge dx_{j_{t+1}}\wedge \cdots \wedge dx_{j_{s}} \\
&\qquad + F f dp\wedge dx_{j_1}\wedge \cdots \wedge dx_{j_s} \\
=&(pdf+\wedge fdp) \wedge\left(
F  dx_{j_1}\wedge\cdots\wedge dx_{j_s} +
 \frac{1}{np}\sum_{t=1}^s x_{j_t}F \wedge dx_{j_1}\wedge\cdots\wedge dx_{j_{t-1}}\wedge dp\wedge dx_{j_{t+1}}\wedge \cdots \wedge dx_{j_{s}} \right) \\
=&dW \wedge \omega'
\end{align*}

We need to show that the map $\bigoplus_{k=0}^\infty(\Omega^s_{\bC[\ud{x}]})_{kn} \to H^0(\bP^{n-1},\pi_*\Omega_{X_{\mathrm{CY}}}^s)$ we defined is bijective. Suppose we have $\omega'\in H^0(\bP^{n-1},\pi_*\Omega_{X_{\mathrm{CY}}}^{r})$ given by the following description:
$$
\omega' := \sum_{1\leq j_1<\cdots<j_s\leq n} F_{(j_1,\ldots,j_s)}(\ud{x},{p})dx_{j_1}\wedge\cdots\wedge dx_{j_s} + \sum_{1\leq j_1<\cdots<j_{s-1}\leq n} G_{(j_1,\ldots,j_{s-1})}(\ud{x},{p})dp\wedge dx_{j_1}\wedge\cdots\wedge dx_{j_{s-1}}.
$$
When the subindex $(j_1,\ldots,j_s)$ is not ordered, we define $\sigma:\{j_1,\cdots,j_s\}\to\{j_1,\cdots,j_s\}$ as the permutation that orders the set, i.e. $\sigma(j_1)<\cdots<\sigma(j_s)$, and define $F_{(j_1,\ldots,j_s)}$ as $F_{(\sigma(j_1),\cdots,\sigma(j_s))}$ times the sign of $\sigma$:
\begin{align}\label{signconvention}
F_{(j_1,\ldots,j_s)} := \mathrm{sign}(\sigma) F_{(\sigma(j_1),\cdots,\sigma(j_s))}.
\end{align}
Rational functioons $F_{(j_1,\ldots,j_s)}$'s and $G_{(j_1,\ldots,j_{s-1})}$'s are related by the fact that $\omega'$ is the kernel of 
$$
\bigwedge^s \left(\{\cO_{X_{\mathrm{CY}}}(-1)\}^{\oplus n}\oplus \cO_{X_{\mathrm{CY}}}(n)\right)\to \Omega_{X_{\mathrm{CY}}}^{s-1} \to 0.
$$
In fact, $G_{(j_1,\ldots,j_{s-1})}$'s are defined by $F_{(j_1,\ldots,j_s)}$'s:
\begin{align}
G_{(j_1,\ldots,j_{s-1})} = \frac{1}{np}\sum_{l=1}^n x_{l} F_{l,j_1,\cdots,j_{s-1}}.
\label{GFrelation}
\end{align}
Recall that $F_{s,j_1,\cdots,j_{s-1}}$ contains the sign from \eqref{signconvention}. If we let $\omega$ as
$$
\omega:=\sum_{1\leq j_1<\cdots<j_s\leq n} F_{j_1,\ldots,j_s}(\ud{x},1)dx_{j_1}\wedge\cdots\wedge dx_{j_s},
$$
then $\omega$ is mapped to $\omega'$ due to the relations \eqref{GFrelation} between $F_{(j_1,\ldots,j_s)}$'s and $G_{(j_1,\ldots,j_{s-1})}$'s.

Thus, we have an isomorphism of complexes between the subcomplex $\left(\bigoplus_{k=0}^\infty\left(\Omega^\bullet_{\bC[\ud{x}]}\right)_{kn},df\right)$ of the Koszul complex and the complex $(H^0(\bP^{n-1},\pi_*\Omega_{X_{\mathrm{CY}}}^\bullet),dW)$.
\end{proof}

\begin{remark}
For the first page in the analytic setting $E_1^{{r},s}=H^s(X_{\mathrm{CY}}^{an},(\Omega^{{r}}_{X_{\mathrm{CY}}})^{an})$, we have a formula similar to \eqref{algcohpage1}:
\begin{align}
H^s(X_{\mathrm{CY}}^{an},(\Omega_{X_{\mathrm{CY}}}^{{r}})^{an}) = \begin{cases}
\bC & \textrm{if $s={r}$, $1\leq {r} \leq n-1,$} \\
H^0(X_{\mathrm{CY}}^{an},(\Omega_{X_{\mathrm{CY}}}^{{r}})^{an}) & \textrm{if $s=0$, ${r}=n$,} \\
0 & \textrm{otherwise}.
\end{cases}
\end{align}
We define $f_1:\widetilde{E}\to E$ as the pullback of $h:(X_{CY}^{an},\cO_{X_{CY}}^{an})\to (X_{CY},\cO_{X_{CY}})$:
\begin{align}
f_1^{{r},s}:H^s(X_{\mathrm{CY}},\Omega_{X_{\mathrm{CY}}}^{{r}}) \to H^s(X_{\mathrm{CY}}^{an},(\Omega_{X_{\mathrm{CY}}}^{{r}})^{an})
\label{spectralsequencemap}
\end{align}
for each ${r},s$.
\end{remark}

\subsection{Algebraic setting versus analytic setting}
\label{sec2.2} 
Here we prove Lemma \ref{lemma5.2}. We start with the definition of an additive grading of $\pi_*\Omega^{{r}}_{X_{\mathrm{CY}}}$, which we will call weight, whose graded pieces are coherent sheaf of $\cO_{\bP_{\bC}^{n-1}}$-modules and $dW$ is homogeneous of weight 1. We define the weight on $\bC[\ud{x},p]$ as follows:
$$
wt(x_1)=\cdots=wt(x_n)=0, \qquad wt(p)=1.
$$
This extends to differential forms on $X_{\mathrm{CY}}$; we define weight of differential forms on $X_{\mathrm{CY}}$ as follows:
\begin{align*}
wt(x_1)=\cdots=wt(x_n)=wt(dx_1)=\cdots=wt(dx_n)=0, \quad
wt(p)=wt(dp)=1.
\end{align*}
More rigorously, from Proposition \ref{lone} and Lemma \ref{lemma5.6}, $\pi_*\Omega^{{r}}_{X_{\mathrm{CY}}}$ is a submodule of
\begin{align*}
\bigwedge^{r} \left(\{\pi_*\cO_{X_{\mathrm{CY}}}(-1)\}^{\oplus n}\oplus \pi_*\cO_{X_{\mathrm{CY}}}(n)\right) \cong  \bigoplus_{k=0}^\infty \bigoplus_{t\in \{0,1\}}\left(\cO_{\bP^{n-1}_{\bC}}\left((k+t)n-{r}+t\right)\right)^{\oplus {n \choose {{r}-t}}}.
\end{align*}
For the component given by fixed $k$ and $t$, we define the \emph{weight} as $k+t$:
\begin{align*}
&\left(\bigwedge^{r} \left(\{\pi_*\cO_{X_{\mathrm{CY}}}(-1)\}^{\oplus n}\oplus \pi_*\cO_{X_{\mathrm{CY}}}(n)\right)\right)_{(l)} \\&:= \left(\cO_{\bP^{n-1}_{\bC}}\left(ln-{r}\right)\right)^{\oplus {n \choose {r}}}\oplus \left(\cO_{\bP^{n-1}_{\bC}}\left((l-1)n-{r}+1\right)\right)^{\oplus {n \choose {{r}-1}}}.
\end{align*}
Then the graded piece $(\pi_*\Omega^{s}_{X_{\mathrm{CY}}})_{(l)}$ of $\pi_*\Omega^{s}_{X_{\mathrm{CY}}}$ is defined as the intersection,
$$
\left(\pi_*\Omega^{{r}}_{X_{\mathrm{CY}}}\right)_{(l)} := \pi_*\Omega^{{r}}_{X_{\mathrm{CY}}} \cap \left(\bigwedge^{r} \left(\{\pi_*\cO_{X_{\mathrm{CY}}}(-1)\}^{\oplus n}\oplus \pi_*\cO_{X_{\mathrm{CY}}}(n)\right)\right)_{(l)}.
$$
From Proposition \ref{lemma5.6}, we have
$$
\pi_* \cO_{X_{CY}} = \bigoplus_{k=0}^\infty \cO_{\bP^{n-1}_{\bC}}(kn),
$$
and this is a sheaf of a graded $\cO_{\bP^{n-1}_{\bC}}$-algebra whose weight $k$ piece is
$$
(\pi_* \cO_{X_{CY}})_{(k)} := \cO_{\bP^{n-1}_{\bC}}(kn).
$$
Then, $\pi_*\Omega^{{r}}_{X_{\mathrm{CY}}}$ is a sheaf of graded $\pi_* \cO_{X_{CY}}$-module, whose multiplication satisfies
$$
m:(\pi_* \cO_{X_{CY}})_{(k)}\otimes_{\cO_{\bP^{n-1}_{\bC}}} (\pi_*\Omega^{{r}}_{X_{\mathrm{CY}}})_{(l)}\to (\pi_*\Omega^{{r}}_{X_{\mathrm{CY}}})_{(k+l)}.
$$
For a fixed $l$, there are only finitely many choices for $k\geq 0$ and $t\in \{0,1\}$ satisfying $k+t=l$. Thus, $(\pi_*\Omega^{{r}}_{X_{\mathrm{CY}}})_{(l)}$ is a subsheaf of a locally free sheaf of finite rank, and so $(\pi_*\Omega^{{r}}_{X_{\mathrm{CY}}})_{(l)}$ is coherent. In addition, $dW$ is a homogeneous operator, since it increases weight by $1$. Since the homogeneous components are all coherent, we can use Serre's GAGA \cite[Th\'eor\`eme 1]{SerreGAGA}. 
Since the cohomology is finite-dimensional, it must be concentrated in the component of weight $\leq l$ for $l$ large enough, which gives the following proposition.
\begin{proposition}
\label{deg<k}
There is a positive integer $l$ such that
$$
H^{r}(H^s(\bP^{n-1}_\bC,\pi_*\Omega_{X_{\mathrm{CY}}}^{\bullet}),dW) \cong H^{r}(H^s(\bP^{n-1}_\bC,(\pi_*\Omega_{X_{\mathrm{CY}}}^{\bullet})_{(\leq l)}),dW).
$$
In particular, for all $k>l$, we have
$$
H^{r}(H^s(\bP^{n-1}_\bC,(\pi_*\Omega_{X_{\mathrm{CY}}}^{\bullet})_{(k)}),dW)=0.
$$
\end{proposition}

\begin{remark}
In the following commutative diagram
\[\begin{tikzcd}
	{X_{\mathrm{CY}}^{an}} & {X_{\mathrm{CY}}} \\
	{\bC\bP^{n-1}} & {\bP_{\bC}^{n-1}}
	\arrow[from=1-1, to=1-2]
	\arrow["{\pi^{an}}", from=1-1, to=2-1]
	\arrow["\pi", from=1-2, to=2-2]
	\arrow["{\iota^{an}}", shift left=3, from=2-1, to=1-1]
	\arrow[from=2-1, to=2-2]
	\arrow["\iota", shift left=3, from=2-2, to=1-2]
\end{tikzcd}\]
where $\iota$ and $\iota^{an}$ are the embeddings along the zero sections, if $\cF$ is a quasi-coherent sheaf of $\cO_{\bP_{\bC}^{n-1}}$-module, $(\iota_* \cF)^{an}$ and $\iota^{an}_* \cF^{an}$ are not the same in general. The sheaf $(\iota_* \cF)^{an}$ consists of functions that are analytic on $X_{\mathrm{CY}}^{an}$, but $\iota^{an}_* \cF^{an}$ consists of functions that are analytic on $\bC\bP^{n-1}$ and algebraic in the direction of the fiber. Thus, in general, we have an injective morphism $\iota^{an}_* \cF^{an}\hookrightarrow (\iota_* \cF)^{an}$ of sheaves of $\cO_{X_{\mathrm{CY}}}^{an}$-modules. Similarly, if $\mathcal{G}$ is a quasi-coherent sheaf of $\cO_{X_{\mathrm{CY}}}$-module, $(\pi^*\mathcal{G})^{an}$ and $\pi^{an}_* \mathcal{G}^{an}$ are not the same in general. There are two different functors $(-)^{an}$: one maps $\mathrm{QCoh}(\cO_{X_{\mathrm{CY}}})$ to $\mathrm{QCoh}(\cO_{X_{\mathrm{CY}}}^{an})$, and the other maps $\mathrm{QCoh}(\cO_{\bP^{n-1}_{\bC}})$ to $\mathrm{QCoh}(\cO_{\bP^{n-1}_{\bC}}^{an})$. We will abusively denote two functors by the same notation $(-)^{an}$. In addition, we will abusively denote $\iota^{an}$ and $\pi^{an}$ by $\iota$ and $\pi$.
\end{remark}

\begin{proposition}
For all $k$,
$$
H^{r}(H^s(\bP^{n-1}_\bC,(\pi_*\Omega_{X_{\mathrm{CY}}}^{\bullet})_{(k)}),dW) \cong H^{r}(H^s(\bC\bP^{n-1},((\pi_*\Omega_{X_{\mathrm{CY}}}^{\bullet})_{(k)}^{an}),dW).
$$
In particular,
$$
H^{r}(H^s(\bP^{n-1}_\bC,(\pi_*\Omega_{X_{\mathrm{CY}}}^{\bullet}),dW) \cong H^{r}(H^s(\bC\bP^{n-1},(\pi_*\Omega_{X_{\mathrm{CY}}}^{\bullet})^{an}),dW).
$$
\label{lemma5.10}
\end{proposition}
\begin{proof}
The first isomorphism on each homogeneous component is derived from Serre's GAGA theorem \cite[Th\'eor\`eme 1]{SerreGAGA} using the fact that $(\pi_*\Omega_{X_{\mathrm{CY}}}^{\bullet})_{(k)}$ is a coherent sheaf on $\bP^{n-1}_{\bC}$. The second isomorphism is obvious, since $\pi_*\Omega_{X_{\mathrm{CY}}}^{\bullet}$ is a direct sum of homogeneous components and the cohomology commutes with the direct sum.
\end{proof}

\begin{proof}[Proof of Lemma \ref{lemma5.2}]
Recall that $\iota:\bC\bP^{n-1}\to X_{\mathrm{CY}}^{an}$ is the closed immersion along the zero section. Then this is an affine morphism, and again by the Leray spectral sequence \cite[Lemma 01F4]{Stacks}, for each ${r},s$, we have
$$
H^{r}(H^s(\bC\bP^{n-1},(\pi_*\Omega_{X_{\mathrm{CY}}}^{\bullet})^{an}),dW) \cong H^{r}(H^s(X_{\mathrm{CY}}^{an},\iota_*(\pi_*\Omega_{X_{\mathrm{CY}}}^{\bullet})^{an}),dW).
$$
From Propositions \ref{leray} and 
\ref{lemma5.10}, for each ${r},s$, we have
$$
H^{r}(H^s(X_{\mathrm{CY}},\Omega^{\bullet}_{X_{\mathrm{CY}}}),dW) \cong H^{r}(H^s(X_{\mathrm{CY}}^{an},\iota_*((\pi_*\Omega_{X_{\mathrm{CY}}}^{\bullet})^{an})),dW).
$$
We want to show that the map $f_1^{{r},s}$ between two cochain complexes 
$$(H^s(X_{\mathrm{CY}}^{an},\iota_*((\pi_*\Omega_{X_{\mathrm{CY}}}^{\bullet})^{an})),dW)\cong(H^s(X_{\mathrm{CY}},\Omega^{\bullet}_{X_{\mathrm{CY}}}),dW)\stackrel{f_1^{\bullet,s}}{\to}(H^s(X_{\mathrm{CY}}^{an},(\Omega_{X_{\mathrm{CY}}}^{\bullet})^{an}),dW)$$
is a quasi-isomorphism; $f_1^{{r},s}$ was defined in \eqref{spectralsequencemap}, and obviously $(dW\wedge)\circ f_1^{{r},s} = f_1^{{r}+1,s}\circ (dW\wedge)$; $\cU=\{U_1,\ldots,U_n\}$ was defined as the open cover of $X_{\mathrm{CY}}$ such that $U_j$ is the inverse image under $\pi$ of the affine open subset in $\bC\bP^{n-1}$ defined by $x_j\neq 0$. In other words, the open set $U_j$ is biholomorphic to $\bC^{n}$ with coordinate functions being $x_1/x_j,\ldots,x_n/x_j,$ $ x_j^{n}p$. Their intersections $U_{i_0,\ldots,i_s}=U_{i_0}\cap \cdots\cap U_{i_s}$ are all Stein. We let $\delta$ be the \v{C}ech differential corresponding to the open cover $\cU$. Thus, by Leray's theorem and Cartan's theorem B \cite[Th\'eor\`eme B]{SerreCartan}, we have the canonical isomorphisms
$$
H^s(X_{\mathrm{CY}}^{an},\iota_*((\pi_*\Omega_{X_{\mathrm{CY}}}^{\bullet})^{an})) \cong \check{H}^s(\cU,\iota_*((\pi_*\Omega_{X_{\mathrm{CY}}}^{\bullet})^{an})),
$$
$$
H^s(X_{\mathrm{CY}}^{an},(\Omega_{X_{\mathrm{CY}}}^{\bullet})^{an}) \cong \check{H}^s(\cU,(\Omega_{X_{\mathrm{CY}}}^{\bullet})^{an}).
$$
Computing $\iota_*((\pi_*\Omega_{X_{\mathrm{CY}}}^{{r}})^{an})$ at $U_{i_0,\ldots,i_s}$, we obtain
\begin{align*}
\iota_*((\pi_*\Omega_{X_{\mathrm{CY}}}^{{r}})^{an})(U_{i_0,\ldots,i_s}) &= (\pi_*\Omega_{X_{\mathrm{CY}}}^{{r}})^{an}(\pi(U_{i_0,\ldots,i_s})) \\
&=\pi_*\Omega_{X_{\mathrm{CY}}}^{{r}} (h\circ\pi(U_{i_0,\ldots,i_s}))\otimes_{\cO_{\bP_{\bC}^{n-1}}(h\circ\pi(U_{i_0,\ldots,i_s}))}\cO_{\bP_{\bC}^{n-1}}^{an}(\pi(U_{i_0,\ldots,i_s}))\\
&=\Omega_{X_{\mathrm{CY}}}^{{r}} (h(U_{i_0,\ldots,i_s}))\otimes_{\cO_{\bP_{\bC}^{n-1}}(h\circ\pi(U_{i_0,\ldots,i_s}))}\cO_{\bP_{\bC}^{n-1}}^{an}(\pi(U_{i_0,\ldots,i_s})).
\end{align*}
Thus, sections of $\iota_*((\pi_*\Omega_{X_{\mathrm{CY}}}^{{r}})^{an})$ are algebraic ${r}$-forms on $X_{\mathrm{CY}}$ tensored with holomorphic functions on $\bC\bP^{n-1}$. In other words, they are holomorphic along the base space $\bC\bP^{n-1}$, and algebraic along the fiber direction. We can interpret them as holomorphic forms on $X_{\mathrm{CY}}^{an}$ via
$$
\Omega_{X_{\mathrm{CY}}}^{{r}} (h(U_{i_0,\ldots,i_s}))\otimes_{\cO_{\bP_{\bC}^{n-1}}(h\circ\pi(U_{i_0,\ldots,i_s}))}\cO_{\bP_{\bC}^{n-1}}^{an}(\pi(U_{i_0,\ldots,i_s})) \to (\Omega_{X_{\mathrm{CY}}}^{{r}})^{an} (U_{i_0,\ldots,i_s}).
$$ 
Thus, we have a map of cochian complexes,
\begin{align}\label{qisom}
(\check{H}^s(\cU,\iota_*((\pi_*\Omega_{X_{\mathrm{CY}}}^{\bullet})^{an})),dW) \to (\check{H}^s(\cU,(\Omega_{X_{\mathrm{CY}}}^{\bullet})^{an}),dW).
\end{align}
An element $[\nu]\in\check{H}^s(\cU,(\Omega_{X_{\mathrm{CY}}}^{p})^{an})$ is represented by a collection of sections
$$
\nu_{i_0,\ldots,i_s}\in (\Omega_{X_{\mathrm{CY}}}^{{r}})^{an}(U_{i_0,\ldots,i_s}),
$$
$$
\nu:=\prod_{i_0,\ldots,i_s}\nu_{i_0,\ldots,i_s}\in \prod_{i_0,\ldots,i_s}(\Omega_{X_{\mathrm{CY}}}^{{r}})^{an}(U_{i_0,\ldots,i_s}),
$$
such that the \v{C}ech differential $\delta\nu$ is zero. We can represent each element $\nu_{i_0,\ldots,i_s}$ as a convergent power series along the fiber direction,
$$
\nu_{i_0,\ldots,i_s} = \sum_{\substack{k\geq 0}} (\eta_{k} P^k+\eta_{k-1}'P^{k-1}\wedge dP)  ,
$$
where $\eta_{k}$ is a holomorphic ${r}$-form on $U_{i_0,\ldots,i_s}$, $\eta_{k}'$ is a holomorphic $({r}-1)$-form on $U_{i_0,\ldots,i_s}$, and $P=x_{i_0}^{n}p$ is the affine coordinate on $U_{i_0,\ldots,i_s}$. 
We define $(\nu_{i_0,\ldots,i_s})_{(k)}$ as the weight $k$ component,
$$
(\nu_{i_0,\ldots,i_s})_{(k)} =\eta_{k} P^k+\eta_{k-1}'P^{k-1}\wedge dP,
$$
$$
\nu_{(k)} = \prod_{i_0,\ldots,i_s}(\nu_{i_0,\ldots,i_s})_{(k)}.
$$

If $dW\wedge \nu=0$, since $dW$ is also homogeneous with respect to the \emph{degree}, we have $dW\wedge \nu_{(k)}=0$. Thus, $\nu_{(k)}$ is representing an element of the cohomology
\begin{align}
\label{qisomissurj}
[\nu_{(k)}]\in H^{r}(H^s(X_{\mathrm{CY}}^{an},\iota_*(\pi_*\Omega_{X_{\mathrm{CY}}}^{\bullet})_{(k)}^{an}),dW).
\end{align}
However, by Proposition \ref{deg<k}, there is a positive integer $l$ such that for $k>l$, 
$$H^{r}(H^s(X_{\mathrm{CY}}^{an},\iota_*(\pi_*\Omega_{X_{\mathrm{CY}}}^{\bullet})_{(k)}^{an}),dW)=0.$$ 
Thus, for $k$ large enough, $\nu_{(k)}= dW\wedge \mu_{(k-1)}$ with $\mu_{(k-1)} \in \check{H}^{s}(\cU,\iota_*(\pi_*\Omega_{X_{\mathrm{CY}}}^{{r}-1})_{(k-1)}^{an})$. Hence, if $dW\wedge \nu=0$, it must be equivalent to the finite sum modulo $dW$,
$$
\nu_{i_0,\ldots,i_s} \equiv \sum_{k=0}^l \eta_{k} P^k+\eta_{k-1}'P^{k-1}\wedge dP \mod dW.
$$
The element in the right hand side sits inside the image of the map \eqref{qisom}. Thus, the map \eqref{qisom} induces a surjection on the cohomologies. 

To show that \eqref{qisom} is injective on the cohomologies, we again observe that $dW$ is weight homogeneous. Suppose the element described in \eqref{qisomissurj} is $0$ in the cohomology $H^{r}(H^s(X_{\mathrm{CY}}^{an},(\Omega^{\bullet}_{X_{\mathrm{CY}}})^{an}),dW)$. In other words, $\nu$ is in the image of $dW$. However, by homogeneity, each graded component $\nu_{(k)}$ must be in the image of $dW$, which implies that $\nu$ represents $0$ in $H^{r}(H^s(X_{\mathrm{CY}},\Omega^{\bullet}_{X_{\mathrm{CY}}}),dW)$. Thus, the kernel of \eqref{qisom} on the cohomology is zero.
\end{proof}

\subsection{Higher terms in the analytic spectral sequence}
\label{sec2.3}

Our next task is to compute the higher differentials $d_k$, $k\geq 2$ in the spectral sequence $E_k^{{r},s}$.

\begin{lemma}
\label{higherdiffan}
The higher differentials $d_k^{{r},s}:E_k^{{r},s}\to E_k^{{r}+k,s-k+1}$, $k>2$ in the spectral sequence are zero.
\end{lemma}
\begin{proof}
When $n$ is even, this is trivial since every higher differential $d_k^{{r},s}$, $k\geq 2$ has either a zero domain or a zero codomain. 

When $n$ is odd and $n\geq 3$, say $n=2{r}+1$ for some ${r}\geq 1$, and the only higher differential that can be possibly nonzero is $d_{{r}+1}^{{r},{r}}: E_{{r}+1}^{{r},{r}}\to E^{n,0}_{{r}+1}$. Since $d_2,\ldots,d_{r}$ are all zeroes, we have
$$
\bC \cong E_{2}^{{r},{r}} \cong E_{{r}+1}^{{r},{r}} \xrightarrow{d_{{r}+1}^{{r},{r}}} E^{n,0}_{{r}+1} \cong E_2^{n,0} \cong \bigoplus_{k=0}^\infty R(f)_{kn}.
$$

We first pick a generator of $E_2^{{r},{r}}$. Let $\omega$ be the Fubini-Study $2$-form which is a generator of $H^{1,1}(\bC\bP^{n-1};\bC)$. Under the identification between the \v{C}ech cohomology and the Dolbeault cohomology of $\bC\bP^{n-1}$, we have a smooth differential form $\theta_j$ on each open subset $U_j\in \cU$ such that
\begin{align}
\theta_j|_{U_j\cap U_k} - \theta_{k}|_{U_j \cap U_k} &= \frac{dx_j}{x_j} - \frac{dx_k}{x_k},
\label{dolbeaultcech}\\
\overline{\partial}\theta_j &= C\cdot\omega,
\label{dolbeaultcech2}
\end{align}
for some nonzero constant $C$. We have the closed immersion $\iota:\bC\bP^{n-1}\to X_{CY}$ which satisfies $\pi\circ \iota$ is the identity on $\bC\bP^{n-1}$. The pullback $(\pi\circ \iota)^*=\iota^*\circ\pi^*$ is also the identity map. Then $\pi^*:H^{{r},{r}}(\bC\bP^{n-1};\bC)\to H^{r}(\cA^{{r},\bullet}(X_{CY}^{an}),\overline{\partial})$ must be an injection, so $\pi^* [\omega^{\wedge {r}}]$ is a generator of $H^{r}(\cA^{{r},\bullet}(X_{CY}^{an}),\overline{\partial})$. Since both $\pi^*\omega^{\wedge {r}}$ and $\omega^{\wedge {r}}$ have the same local description, we abuse the notation $\omega^{\wedge {r}}$ for $\pi^*\omega^{\wedge {r}}$.

Dolbeault's isomorphism $H^s(\cA^{{r},\bullet}(X_{\mathrm{CY}}^{an}),\overline{\partial})\cong H^{s}(X_{\mathrm{CY}}^{an},(\Omega^{{r}}_{X_{\mathrm{CY}}})^{an})$ maps $C\cdot\omega$ to
$$
\left[\left(\frac{dx_j}{x_j} - \frac{dx_k}{x_k}\right)_{(j,k)}\right]\in H^{1}(X_{\mathrm{CY}}^{an},(\Omega_{X_{\mathrm{CY}}}^{1})^{an})\cong H^{1}(X_{\mathrm{CY}},\Omega_{X_{\mathrm{CY}}}^{1}).
$$
To be specific, we define $\theta_j$ as follows:
$$
\theta_j := \partial \log\left(1+\sum_{\substack{1\leq l\leq n\\l\neq j}}|X_{l/j}|^2\right)=\sum_{\substack{1\leq k\leq n\\k\neq j}}\frac{\overline{X}_{k/j}dX_{k/j}}{1+\sum_{\substack{1\leq l\leq n\\l\neq j}}|X_{l/j}|^2},
$$
where $X_{k/j}= x_k/x_j$ is the affine coordinate of $U_j$.

A crucial equality for the proof of $d_{{r}+1}^{{r},{r}}=0$ is the following: for $j,k\in \{1,\ldots,n\}$,
\begin{align}
\frac{1}{n}\left(d\left(\frac{f}{x_j^{n}}\right)\wedge d(x_j^{n}p)-d\left(\frac{f}{x_k^{n}}\right)\wedge d(x_k^{n}p)\right) = dW\wedge \left(\frac{dx_j}{x_j}-\frac{dx_k}{x_k}\right),
\label{higherdiffeq}
\end{align}
which is derived from the following elementary equation: for each $j=1,\cdots,n$,
\begin{align}
\frac{1}{n} d\left(\frac{f}{x_j^{n}}\right)\wedge d(x_j^{n}p) = \frac{1}{n} df\wedge dp+ dW \wedge \frac{dx_j}{x_j}.
\label{higherdiffeq2}
\end{align}
An element $(\frac{dx_j}{x_j}-\frac{dx_k}{x_k})_{j,k}\in \check{C}^1(\cU,\Omega_{X_{\mathrm{CY}}}^{1})$ is a generator of $H^1(X_{\mathrm{CY}},\Omega_{X_{\mathrm{CY}}}^{1})\cong H^1(X_{\mathrm{CY}}^{an},(\Omega_{X_{\mathrm{CY}}}^{1})^{an})$. We have an element $\left(\frac{1}{n}{d\left(\frac{f}{x_j^{n}}\right)\wedge d(x_j^{n}p)}\right)_{j}$ of $\check{C}^0(\cU,\Omega_{X_{\mathrm{CY}}}^{2})$, which is not in the kernel of $\delta$. However, after taking $dW$ we get zero:
\begin{align}
\label{higherdiffeq3}
dW\wedge \frac{1}{n}{d\left(\frac{f}{x_j^{n}}\right)\wedge d(x_j^{n}p)} = 0.
\end{align}
We interpret the equations \eqref{higherdiffeq} and \eqref{higherdiffeq2} in the spectral sequence $E_k^{{r},s}$ to show that $d_{{r}+1}^{{r},{r}}=0$. 

At $U_j$, we consider the $(2,0)$-form $\alpha$ whose restriction at $U_j$ is
\begin{align}
\alpha|_{U_j}=\frac{1}{n}d\left(\frac{f}{x_j^n}\right)\wedge d(x_j^{n}p) - dW \wedge \theta_j.
\label{defalpha}
\end{align}
This is well-defined global $(2,0)$-form since by \eqref{higherdiffeq} and \eqref{dolbeaultcech},
$$
\alpha|_{U_j \cap U_k} - \alpha|_{U_k \cap U_j} = 0.
$$
Moreover, by \eqref{higherdiffeq3} and \eqref{dolbeaultcech2}, we have
\begin{align}
\label{alphaprop1} dW\wedge \alpha = 0, \\
\label{alphaprop2} \overline{\partial}\alpha = CdW\wedge \omega.
\end{align}
Thus, along the stair case
\[\begin{tikzcd}
	{\cA^{1,1}(X_{\mathrm{CY}}^{an})} & {\cA^{2,1}(X_{\mathrm{CY}}^{an})} \\
	& {\cA^{2,0}(X_{\mathrm{CY}}^{an})} & {\cA^{3,0}(X_{\mathrm{CY}}^{an}),}
	\arrow["dW", from=1-1, to=1-2]
	\arrow["{\overline{\partial}}"', from=2-2, to=1-2]
	\arrow["dW", from=2-2, to=2-3]
\end{tikzcd}\]
the Fubini-Study $2$-form $\omega$ is transferred as follows:
\[\begin{tikzcd}
	{C\cdot\omega} & {CdW\wedge \omega} \\
	& \alpha & {0.}
	\arrow["dW", maps to, from=1-1, to=1-2]
	\arrow["{\overline{\partial}}"', maps to, from=2-2, to=1-2]
	\arrow["dW", maps to, from=2-2, to=2-3]
\end{tikzcd}\]
Applying this relation to $C^{r} \omega^{\wedge {r}}$ which is a generator of $E_{{r}+1}^{{r},{r}}$, we have
\[\begin{tikzcd}
	{C^{r} \omega^{\wedge {r}}} & {C^{{r}}dW\wedge \omega^{\wedge {r}}} \\
	& C^{{r}-1}\alpha \wedge \omega^{\wedge ({r}-1)}& {0.}
	\arrow["dW", maps to, from=1-1, to=1-2]
	\arrow["{\overline{\partial}}"', maps to, from=2-2, to=1-2]
	\arrow["dW", maps to, from=2-2, to=2-3]
\end{tikzcd}\]
Hence,
$$
d_{{r}+1}^{{r},{r}}(C^{r} \omega^{\wedge {r}}) = 0.
$$
\end{proof}

\section{Comparison of Frobenius algebra structures}\label{sec3}
We introduce two dGBV algebras and two Frobenius algebras. We provide a precise comparison result between them: we prove Theorem \ref{theoremtwo}.

\subsection{Two Frobenius algebras: Barannikov-Kontsevich versus Li-Wen}\label{sec3.1}

We start from recalling two Frobenius algebra structures on $H(V;\C)$ due to \cite{BK} and \cite{LW}.
By the Hodge decomposition $H^k(V;\C) \cong \bigoplus_{p+q=k} H^{p,q}(V)$ where $H^{p,q}(V)= H^q(V,\wedge^p\cT_V^*)$.
The contraction with a non-vanishing holomorphic volume $(n-2)$-form $\O_V$ on $V$ induces an isomorphism
\be
 \U_{BK}: \bigoplus_{i,j=0}^{n-2} \cA^{0,j}(V, \wedge^i \cT_V) \xrightarrow{\cong} \bigoplus_{i,j=0}^{n-2}\cA^{0,j}(V, \wedge^{n-2-i}\cT^*_V)
\ee
and hence an isomorphism $H^j(V, \wedge^i \cT_V )\cong H^{n-2-i,j}(V)$ for each $0\leq i,j\leq n-2$ after taking $\overline{\partial }$-cohomology. 
Note that the space $ \bigoplus_{i,j=0}^{n-2} \PV^{i,j}(V)$ has a product structure coming from the product of holomorphic polyvector fields valued in anti-holomorphic differential forms:
\be \label{CYR}
\wedge \colon  \cA^{0,j}(V, \wedge^i \cT_V)\times \cA^{0,l}(V, \wedge^k \cT_V) \to \cA^{0,j+l}(V, \wedge^{i+k} \cT_V).
\ee
Recall that
\be \label{pvgrading}
\PV^r(V):= \bigoplus_{r=j-i}\PV^{i,j}(V).
\ee
Then in the case of $r=0$, the contraction map $ c_{\O_V}$ restricts to an isomorphism
\be
 \U_{BK}: \bigoplus_{j=0}^{n-2} \cA^{0,j}(V, \wedge^{j} \cT_V) \xrightarrow{\cong} \bigoplus_{j=0}^{n-2}\cA^{0,j}(V, \wedge^{n-j}\cT^{*}_V)
\ee
and induces an isomorphism $ H^0(\PV(V) ):=\bigoplus_{j=0}^{n-2}  H^{j,j}(\PV(V))\xrightarrow{\cong} H^{n-2}(V)$ where $ H^{j,j}(\PV(V)):=H^j(V, \wedge^j \cT_V)$. Under this isomorphism, $H^{j,j}(\PV(V))$ corresponds to $H^{n-2-j,j}(V)$.
In addition, it has a trace map $\Tr_{BK}^{\O_V}: \PV(V) \to  \C$ given by
\bea\label{tracecy}
\Tr_{BK}^{\O_V}(\o) = \int_V \left(\o \vdash \O_V\right) \wedge \O_V.
\eea
Clearly, the trace map vanishes unless $\omega \in \PV^{n-2,n-2}(V)$ and, in fact, it induces an isomorphism
$ H^{n-2,n-2}(\PV(V))  \cong \C$. 
Using $\Tr_{BK}^{\O_V}$, we define a symmetric bilinear pairing $K_{BK}^{(0)}$ on $H^0(\PV(V),\bar\partial)$ by
\bea \label{cypair}
K_{BK}^{(0)}(\omega,\eta):=\Tr_{BK}^{\O_V}(w\wedge  \eta)= \int_V \left((\o \wedge  \eta ) \vdash \O_V\right) \wedge \O_V,
\eea
for $\o, \eta \in \PV(V)$. This yields a Frobenius algebra $( \bigoplus_{i,j=0}^{n-2} H^j(V , \wedge^i \cT_V )  , \wedge ,  K_{BK}^{(0)})$.

%

\begin{definition}\label{BKF}
The BK Frobenius algebra $(\O_{BK}, K_{BK}^{(0)})$ is defined to be
\be
( H(\PV(V),\bar\partial),\wedge, K_{BK}^{(0)} )=( \bigoplus_{i,j=0}^{n-2} H^j(V , \wedge^i \cT_V ) , \wedge , K_{BK}^{(0)}).
\ee
\end{definition}

Now we review the LW Frobenius algebra due to \cite{LW}.
We have natural inclusions of complexes
\[\begin{tikzcd}
	{(\PV_c(X_{\mathrm{CY}}),\bar\partial_W)} & {(\PV_{W,\infty}(X_{\mathrm{CY}}),\bar\partial_W)} & {(\PV(X_{\mathrm{CY}}),\bar\partial_W)}
	\arrow["{\iota_1}", hook, from=1-1, to=1-2]
	\arrow["{\iota_2}", hook, from=1-2, to=1-3]
\end{tikzcd}\]
where $\PV_{W,\infty}(X_{\mathrm{CY}})$ is the $W$-twisted Sobolev space of polyvector fields
(\cite[Definition 3.4]{LW}) and $\PV_c(X_{\mathrm{CY}})$ is the space of polyvector fields with compact support.
The space $ \bigoplus_{i,j=0}^n \PV^{i,j}(X_{\mathrm{CY}})$ also has a product structure coming from the product of holomorphic polyvector fields valued in anti-holomorphic differential forms:
\be 
\wedge \colon  \cA^{0,j}(X_{\mathrm{CY}}, \wedge^i \cT_{X_{\mathrm{CY}}})\times \cA^{0,l}(X_{\mathrm{CY}}, \wedge^k \cT_{X_{\mathrm{CY}}}) \to \cA^{0,j+l}(X_{\mathrm{CY}}, \wedge^{i+k} \cT_{X_{\mathrm{CY}}}).
\ee
The same product structure is induced on $ \PV_{W,\infty}(X_{\mathrm{CY}})=\bigoplus_{i,j=0}^n \PV^{i,j}_{W,\infty}(X_{\mathrm{CY}})$.

 Since $(X_{\mathrm{CY}},g, \Omega_{\mathrm{CY}})$ is a bounded Calabi-Yau geometry and $W$ satisfies the strongly elliptic condition, both $\iota_1$ and $\iota_2$ are quasi-isomorphisms
due to Theorem \cite[Theorem 3.14]{LW}. Let $\rho$ be a smooth function with compact support such that $\rho =1$ in a neighborhood of $\Crit(W)$. The cohomological inverse can be described by using two operators $V_W:\cA(X_{CY}\setminus \Crit(W))\to \cA(X_{CY}\setminus \Crit(W))$ and $T_\rho:\cA(X_{CY})\to \cA_{W,\infty}(X_{CY})$. They are defined as follows:
\begin{align}
V_W &:= \sum_{j,k} \frac{\overline{W_j}g^{\overline{j},k}}{\|\nabla W\|^2}\iota_{\partial_{k}}: \quad \cA^{*,*}(X_{CY}\setminus \Crit(W)) \to \cA^{*-1,*}(X_{CY}\setminus \Crit(W)) \, \label{defVW}\\
T_\rho &:=\rho + [\overline\partial, \rho] V_W \sum_{j=0}^\infty (-1)^j[\overline{\partial},V_W]^j= \rho + (\overline\partial \rho) V_W \sum_{j=0}^\infty (-1)^j[\overline{\partial},V_W]^j, \label{defTrho}
\end{align}
which is well-defined since $[\overline{\partial},V_W]^j=0$ for $j$ large enough. Here $\iota_{\partial_{k}}$ represents the contraction with respect to $\partial_{k}$, and $(g^{\overline{j},k})$ is the inverse matrix of the matrix $(g_{j,\overline{k}})$ representing the K\"ahler metric on $X_{CY}$:
$$
g_{j,\overline{k}} = g(\partial_j,\overline{\partial_k}).
$$
Then, according to the proof of Theorem \cite[Theorem 3.14]{LW}, $T_{\rho}$ is a quasi-isomorphism from $(\cA(X_{CY}),\overline{\partial}_W)$ to $(\cA_{W,\infty}(X_{CY}),\overline{\partial}_W)$, which is the cohomological inverse to the embedding $\iota_2$: from now own, we identify $H(\PV_{W,\infty}(X_{CY}), \bar\partial_W)$ with $H(\PV(X_{CY}), \bar\partial_W)$ by using $\iota_2$ and $T_\rho$.
According to \cite[Lemma 3.7]{LW}, the contraction map with $\O_{\mathrm{CY}}$ induces an isomorphism $\U_{LW}:\PV_{W,\infty}(X_{\mathrm{CY}}) \to \cA_{W,\infty}(X_{\mathrm{CY}})$, where $\cA_{W,\infty}(X_{\mathrm{CY}})$ is the space of $W$-twisted Sobolev differential forms (\cite[Definition 2.19]{LW}).
Moreover, \cite[Theorem 3.9]{LW} says that $(\PV_{W,\infty}(X_{\mathrm{CY}}), \bar\partial_W, \O_{\mathrm{CY}})$ also forms a dGBV algebra.

Following \cite[Definition 3.16]{LW}, we define 
\be
K_{LW} :\PV_{W,\infty}(X_{\mathrm{CY}})[[u]] \times \PV_{W,\infty}(X_{\mathrm{CY}})[[u]]  \to \C[[u]]
\ee
by
\be
K_{LW}(a(u) \a, b(u) \b):= a(u)b(-u) \int_{X_{\mathrm{CY}}} \U_{LW}(\a\wedge\b) \wedge \O_{\mathrm{CY}}.
\ee

Note that $\bar\partial_W$ is graded skew-symmetric and $\partial$ is graded symmetric with respect to the pairing $K_W$:
\be
K_{LW}(\bar\partial_W \a, \b) &=& (-1)^{|\a|+1} K_{LW}(\a, \bar\partial_W \b),\\
K_{LW}(\partial_{\mathrm{CY}} \a, \b) &=& (-1)^{|\a|} K_{LW}(\a, \partial_{\mathrm{CY}} \b),
\ee
for all $\a, \b \in \PV_{W,\infty}(X_{\mathrm{CY}}).$
Thus $K_{LW}$ induces a sesquilinear pairing \cite[Definition 3.18]{LW}
\be
K_{LW}:H(\PV_{W,\infty}(X_{\mathrm{CY}})[[u]], \bar\partial_W+u\partial_{\O_{\mathrm{CY}}})\times H(\PV_{W,\infty}(X_{\mathrm{CY}})[[u]], \bar\partial_W+u\partial_{\O_{\mathrm{CY}}}) \to \C[[u]]
\ee
and $K_{LW}=\sum_{k\geq 0} u^k K_{LW}^{(k)}$ induces a non-degenerate $\C$-bilinear pairing \cite[Theorem 3.20]{LW}
\be
K_{LW}^{(0)}: H(\PV_{W,\infty}(X_{\mathrm{CY}}), \bar\partial_W) \times H(\PV_{W,\infty}(X_{\mathrm{CY}}), \bar\partial_W)
\to \C.
\ee

Then the triple $(H(\PV_{W,\infty}(X_{\mathrm{CY}}),\bar\partial_W), \wedge ,  K_{LW}^{(0)})$ yields a Frobenius algebra.
\begin{definition}\label{LWF}
The LW Frobenius algebra $(\O_{LW}, K_{LW}^{(0)})$ is defined to be
\be
(H(\PV_{W,\infty}(X_{\mathrm{CY}}),\bar\partial_W), \wedge , K_{LW}^{(0)}).
\ee
\end{definition}

Note that
$$
K_{LW}^{(0)}([\U_{LW}^{-1}\circ T_\rho(\alpha)],[\U_{LW}^{-1}\circ T_\rho(\beta)])=\int_{X_{CY}}\U_{LW}(\U_{LW}^{-1}\circ T_\rho(\alpha)\wedge \U_{LW}^{-1}\circ T_{\rho}(\beta))\wedge \Omega_{CY},
$$
for $\a, \b \in \PV(X_{CY})$. Similar to the BK side, we define $\Tr_{LW}:PV_{W,\infty}(X_{CY})\to \bC$ as
$$
\Tr_{LW}(\eta) := K^{(0)}_{LW}(1,\eta) = \int_{X_{CY}}\U_{LW}(\eta)\wedge\Omega_{CY}
$$
for $\eta\in PV_{W,\infty}(X_{CY})$, which induces a well-defined map $$
\widehat{\Tr}_{LW}:H(PV_{W,\infty}(X_{CY}),\overline{\partial}_W)=H(PV(X_{CY}),\overline{\partial}_W)\to \bC.
$$


\subsection{Comparison theorem}\label{sec3.2}
We prove Theorem \ref{theoremtwo}, i.e. we will construct a ring isomorphism
    \begin{eqnarray*}
        \Phi:\O_{LW}=(H(\PV_{W,\infty}(X_{\mathrm{CY}}), \bar\partial_W), \wedge) \to \O_{BK}=(H(\PV(V), \bar\partial), \wedge)
    \end{eqnarray*}
    such that there exists a non-zero constant $c_\Phi \in \C$ satisfying
    \begin{eqnarray}\label{ip}
        c_\Phi K_{LW}^{(0)}(\a, \b) = K_{BK}^{(0)}(\Phi(\a), \Phi(\b)), \quad \a, \b \in \O_{LW}.
    \end{eqnarray}

From Theorem \ref{cohomologycmp} and the computation of the spectral sequence in \eqref{secondpage}, we obtain an injective ring homomorphism
$$
\Ker\overline{\partial}\cap PV^{0,0}(X_{CY})\supseteq \bC[\ud x,p]_0\twoheadrightarrow R(W)_0 \to H^0(PV(X_{CY}),\overline{\partial}_W),
$$
and the short exact sequence of $\bC$-vector spaces
$$
0\to R(W)_0 \to H(PV(X_{CY}),\overline{\partial}_W) \to \bC^{n-1} \to 0.
$$
 For simplicity, write $\omega_0=\overline{\partial}\theta_j$; see \eqref{dolbeaultcech2}. Then, for $q=1,\ldots,n-1$, the differential form $\omega_0^{\wedge q}-\alpha\wedge\omega_0^{\wedge(q-1)}$ is an element of $\cA^{2q}(X_{CY})$; see \eqref{defalpha} for the definition of $\alpha$. Let us define 
\begin{align}
\alpha_{2q-n} := \U_{LW}^{-1}(\omega_0^{\wedge q}-\alpha\wedge\omega_0^{\wedge(q-1)})\in PV^{n-q,q}(X_{CY}).
\label{nonprgenerator}
\end{align}
This is a representative of a generator of $H^{2q-n}(PV(X_{CY}),\overline{\partial}_W)$ by \eqref{alphaprop1} and \eqref{alphaprop2}.
We also observe that for $j=0,\ldots,n-2$ and $g(\ud{x})\in \bC[\ud{x}]_{jn}$, $[p^jg(\ud{x})]$ is an element in $H^{0}(PV(X_{CY}),\overline{\partial}_W)$.

Let $e_1, \cdots, e_{n-1}$ to be the standard basis of $\C^{n-1}$. We describe a map $\Psi_{LW}$ on $\C^{n-1}$:
\begin{eqnarray}\label{dpsi}
    \Psi_{LW}:\C^{n-1} \to H(\PV(X_{CY}), \bar\partial_W), \quad e_i \to \Psi_{LW}(e_i)=[\a_{2i-n}], \ i=1, \ldots, n-1.
\end{eqnarray}

\begin{definition}
We define a vector space isomorphism $\Psi_{LW}:R(W)_0\oplus\bC^{n-1}\to H(PV(X_{CY}),\overline{\partial}_W)$ as an extension of \eqref{dpsi} via the ring monomorphism $R(W)_0 \to H(PV(X_{CY}),\overline{\partial}_W)$.
\end{definition}

Then for $[p^jg(x)]\in R(W)_0$, obviously, $\Psi_{LW}([p^jg(\ud{x})])=[p^jg(\ud{x})]\in H^0(PV(X_{CY}),\overline{\partial}_W)$. 
We have the following lemma on the multiplication structure of $H(PV(X_{CY}),\overline{\partial}_W)$.
\begin{lemma}\label{multiplication1}
For $j>0$ and $p^jg(\ud{x})\in \cO_{X_{CY}}(X_{CY})$, we have $[p^jg(\ud{x})\alpha_q]=0$ for all $q=-n+2,-n+4,\ldots,n-2$.
\end{lemma}
\begin{proof}
Since $\U_{LW}(p^{j}g(\ud{x})\alpha_q)=p^jg(\ud{x})\U_{LW}(\alpha_q)$, it is enough to show that for each $q=1,\ldots,n-1$ and $j>0$, $p^jg(\ud{x})(\omega_0^{\wedge q}-\alpha\wedge\omega_0^{\wedge(q-1)})$ is zero in $H^{2q}(\cA(X_{CY}),\overline{\partial}+dW)$. Locally on the open set $U_k$, we have
$$
(\overline{\partial}+dW)\left(P^jg(\ud{X})\theta_k\wedge\omega_0^{\wedge(q-1)}-\frac{1}{n}P^{j-1}g(\ud{X})dP\wedge\omega_{0}^{\wedge(q-1)}\right) = P^jg(\ud{X})(\omega_0^{\wedge q}-\alpha\wedge\omega_0^{\wedge(q-1)}).
$$
We only need to prove that the collection of differential forms $\beta_k$ defined on each open set $U_k$,
$$
\beta_k = P^jg(\ud{X})\theta_k\wedge\omega_0^{\wedge (q-1)}-\frac{1}{n}P^{j-1}g(\ud{X})dP\wedge\omega_{0}^{\wedge(q-1)}
$$
defines a global differential form. Let $(Y_1,\ldots,Y_{l-1},Y_{l+1},\ldots,Y_n,Q)$ be the coordinate function on the open set $U_l$, so that
\begin{align*}
X_j/X_l = Y_j,\ (j\neq k,l)& & 1/X_l=Y_k, && X_l^nP=Q. 
\end{align*}
Then, using the above transition function,
$$
g(\ud{X})P^{j-1}dP-g(\ud{Y})Q^{j-1}dQ = nP^jg(\ud{X})\frac{dX_l}{X_l}=nP^jg(\ud{X})\left(\frac{dx_l}{x_l}-\frac{dx_k}{x_k}\right).
$$
From \eqref{dolbeaultcech}, we obtain
$$
\beta_k|_{U_k\cap U_l}- \beta_l|_{U_k\cap U_l}=0.
$$
Thus, $(\beta_k)$ defines a global differential form $\beta\in \cA^{2q-1}(X_{CY})$ such that $(\overline{\partial}+dW)\beta = p^jg(\ud{x})(\omega_0^{\wedge q}-\alpha\wedge\omega_0 ^{\wedge(q-1)})$.
\end{proof}

\begin{lemma} \label{integrationcompute}
The trace $\Tr_{LW}$ on $H(PV_{W,\infty}(X_{CY}),\overline{\partial}_W)$ has the following properties:
\begin{enumerate}
\item
For $0\leq j\leq n-3$ and $g(\ud{x})\in \bC[\ud{x}]_{nj}$, we have $\Tr_{LW}(T_{\rho}(p^jg(\ud{x})))=0$.
\item
For $q=-n+2,-n+4,\ldots,n-2$, $\Tr_{LW}(T_{\rho}(\alpha_q))=0$.
\item
There exists a homogeneous polynomial $h(\ud{x})\in \bC[\ud{X}]_{n(n-2)}$ such that $\Tr_{LW}(T_\rho(p^{n-2}h(\ud{x})))\neq 0$.
\end{enumerate}
\end{lemma}
\begin{proof}
The results follow from the direct computation of $\Tr_{LW}(T_\rho(p^j g(\ud{x})))$ and $\Tr_{LW}(T_\rho(\alpha_q))$.

(1) Recall the definition of $V_W$ and $T_\rho$ given in \eqref{defVW} and \eqref{defTrho}. First we compute the integral,
$$
\Tr_{LW}(T_\rho(p^{j}g(\ud{x})))=(-1)^{n-1}\int_{X_{CY}} (\overline{\partial}\rho)(V_W [\overline{\partial},V_W]^{n-1}(T_\rho(p^{j}g(\ud{x})\Omega_{CY}))\wedge \Omega_{CY}.
$$
Since $p^j g(\ud{x})$ is holomorphic,
\begin{align*}
\Tr_{LW}(T_\rho(p^{j}g(\ud{x})))&=(-1)^{n-1}\int_{X_{CY}} p^{j}g(\ud{x})(\overline{\partial}\rho)\left(V_W [\overline{\partial},V_W]^{n-1}( \Omega_{CY})\right)\wedge \Omega_{CY}\\
&=(-1)^{n-1}\int_{U_n} P^{j}g(\ud{X})(\overline{\partial}\rho)\left(V_W [\overline{\partial},V_W]^{n-1}(\Omega_{CY})\right)\wedge \Omega_{CY}.
\end{align*}
Let $A:=V_W [\overline{\partial},V_W]^{n-1}(\Omega_{CY})$, $X':=$ Support of $\overline{\partial}\rho$. Then,
\begin{align*}
\int_{U_n} P^j g(\ud{X})(\overline{\partial}\rho)\left(V_W [\overline{\partial},V_W]^{n-1}(\Omega_{CY})\right)\wedge \Omega_{CY} &=  \int_{X'\cap U_n} P^j g(\ud{X}) (\overline{\partial}\rho)A\wedge \Omega_{CY}
\end{align*}
Let $P=re^{i\theta}$ for $r\in \bR$, $\theta\in [0,2\pi)$. Since the metric on $X_{CY}$ is given so that $X_{CY}\to \bC^n/\mu_n$ is a crepent resolution, we can take $\Omega_{CY}=dX_1\wedge\cdots\wedge dX_{n-1}\wedge dP$. It is possible to pick $\rho$ so that ${\partial}\rho/\partial\theta=0$; consider the following radius function
$$
R^{2n} = {|p|^2}{(|x_1|^2+\cdots+|x_n|^2)^n}\geq 0
$$
defined on $X_{CY}$. Inequality $R\leq 1$ defines a compact set in $X$. We define
$$
\rho := \begin{cases}
1 & \textrm{if $R\leq 1,$}\\
\frac{1}{1+e^{-\frac{1}{R^{2n}-1}+\frac{1}{2-R^{2n}}}}& \textrm{if $1<R<2^{\frac{1}{2n}}$,} \\
0 & \textrm{if $R\geq 2^{\frac{1}{2n}}$.}
\end{cases}
$$ Then,
$$
\overline{\partial}\rho = \sum_{j=1}^{n-1}\frac{\partial \rho}{\partial\overline{X}_j}d\overline{X}_j + \frac{\partial \rho}{\partial r}\frac{e^{i\theta}}{2} d\overline{P}=\sum_{j=1}^{n-1}\frac{\partial \rho}{\partial\overline{X}_j}d\overline{X}_j + \frac{1}{2}\frac{\partial \rho}{\partial r}(dr-ird\theta),
$$
\begin{align*}
\int_{X'\cap U_n} P^j g(\ud{X})(\overline{\partial}\rho) A\wedge \Omega_{CY} = \int_{\rho(X_1,\ldots,X_{n-1},1,r)\neq 0}\int_{0\leq \theta \leq 2\pi} P^j g(\ud{X}) (\overline{\partial}\rho)A\wedge \Omega_{CY}.
\end{align*}

 We let $Z_0=P$, $Z_j=X_j$ for $j=1,\ldots,n-1$ and $W_k=\frac{\partial}{\partial Z_k} W$. Similar to the computation in Li-Li-Saito \cite[Proposition 2.5]{LLS}, we have the following
\begin{align*}
A&= \frac{C}{\|\nabla W\|^{2n}}\left(\sum_{ j_0\neq\cdots\neq j_{n-1}}\sum_{k=0}^{n-1}(-1)^{k}\left(\overline{W_k}g^{\overline{k},j_k}\right) \overline{\partial}\left(\overline{W_0}g^{\overline{0},j_0}\right)\wedge \cdots \wedge \widehat{\overline{\partial}\left(\overline{W_k}g^{\overline{k},j_k}\right)} \wedge \cdots \wedge \overline{\partial}\left(\overline{W_{n-1}}g^{\overline{n-1},j_{n-1}}\right)\right),
\end{align*}
for some constant $C$.  Our next goal is to compute the order of $\overline{P}$ and $P$ in $g^{j\overline{k}}$.

From the work of D. Joyce \cite[Example 3.5]{Joyce}, the metric is given by the following K\"ahler potential:
\begin{align*}
F&= \sqrt[n]{R^{2n}+1}+\frac{1}{n}\sum_{j=0}^{n-1} \zeta^j \log\left(\sqrt[n]{R^{2n}+1}-\zeta^j\right),\\
\eta&=i\partial\overline{\partial}F.
\end{align*}
On the open set $U_n$, $\partial F/\partial \theta=0$. Then,
$$
g_{{j},\overline{k}} = \frac{\partial^2 F}{\partial {Z}_j\partial \overline{Z}_k},
$$
and $g^{\overline{j},k}$ is a $(j,k)$-entry of the inverse matrix of $(g_{{j},\overline{k}})_{j,k}$. If $j\neq 0$,
\begin{align*}
g_{j,\overline{0}} &= P\frac{1}{2r}\frac{\partial^2F}{\partial X_j\partial r},\\
g_{0,\overline{j}} &= \overline{P}\frac{1}{2r}\frac{\partial^2F}{\partial \overline{X}_j\partial r},\\
g_{0,\overline{0}} &= \frac{r}{4}\frac{\partial^2 F}{\partial r^2} + \frac{1}{4r} \frac{\partial F}{\partial r}.
\end{align*}

The determinant of $(g_{{j},\overline{k}})$ is constant by the fact that $\omega^{\wedge n} = c\cdot\Omega\wedge \overline{\Omega}$ for some constant $c$. To compute $g^{\overline{j},k}$, we observe that the matrix $(g_{{j},\overline{k}})$ is in the following form
$$
\begin{bmatrix}
h_{0,0} & \overline{P}h_{0,1} & \cdots & \overline{P}h_{0,n-1} \\
{P} h_{1,0} & h_{1,1} & \cdots & h_{1,n-1}  \\
\vdots & \vdots &\ddots & \vdots \\
{P} h_{n-1,0} & h_{n-1,1} & \cdots & h_{n-1,n-1}
\end{bmatrix},
$$
where $h_{i,j}$ is functions irrelevant to $\theta$.
Since $(g^{\overline{j},k})$ is the inverse matrix of $(g_{{j},\overline{k}})$, the value $g^{\overline{j},k}$ is given by the determinant of the cofactor matrix $M_{k,j}$ of $(g_{{j},\overline{k}})$ times $(-1)^{j+k}$.

If $(j,k)=(0,0)$, we have
$$
g^{\overline{0},0} = \frac{1}{\det (g_{j,\overline{k}})}\det M_{0,0}.
$$
The determinant of a cofactor matrix $M_{0,0}$ is irrelevant to $\theta$, since it is an algebraic function over $h_{j,k}$'s with $j,k\neq 0$. Thus, $g^{{0},\overline{0}}$ is irrelevant to $\theta$. 

When $j=0$ and $k\neq 0$, 
$$
g^{\overline{0},{k}} =  \frac{(-1)^k}{\det (g_{j,\overline{k}})}\det M_{k,0}.
$$
Then the first row of $M_{k,0}$ is
$$
\begin{bmatrix}
\overline{P}h_{0,1} &\cdots & \overline{P} h_{0,n-1}
\end{bmatrix}.
$$
By the property of the determinant, we have
$$
g^{\overline{0},k} = \frac{(-1)^k}{\det (g_{j,\overline{k}})} \overline{P}\det\begin{bmatrix}
h_{0,1} & \cdots & h_{0,n-1} \\
\vdots & \ddots &\vdots   \\
h_{k-1,1} & \cdots & h_{k-1,n-1} \\
h_{k+1,1} & \cdots & h_{k+1,n-1} \\
\vdots & \ddots &\vdots   \\
h_{n-1,1} & \cdots & h_{n-1,n-1} 
\end{bmatrix}.
$$
The determinant of the matrix above is irrelevant to $\theta$, so we conclude that $g^{\overline{0},k}$ is $\overline{P}$ times a function $h^{0,k}$ which is irrelevant to $\theta$. Similarly, $g^{\overline{j},0}$ is ${P}$ times a function $h^{j,0}$ which is irrelevant to $\theta$.

When $j\neq 0$ and $k\neq 0$, the determinant of $M_{k,j}$ is
\begin{align*}
\det\begin{bmatrix}
h_{0,0} & \overline{P}h_{0,1} & \cdots & \overline{P}h_{0,j-1} & \overline{P}h_{0,j+1} & \cdots& \overline{P}h_{0,n-1} \\
{P} h_{1,0} & h_{1,1} & \cdots & h_{1,j-1} & h_{1,j+1}& \cdots & h_{1,n-1}  \\
\vdots & \vdots &\ddots & \vdots & \cdots &\ddots  & \vdots \\
{P} h_{k-1,0} & h_{k-1,1} & \cdots & h_{k-1,j-1} & h_{k-1,j+1}& \cdots & h_{k-1,n-1}  \\
{P} h_{k+1,0} & h_{k+1,1} & \cdots & h_{k+1,j-1} & h_{k+1,j+1}& \cdots & h_{j+1,n-1}  \\
\vdots & \vdots &\ddots & \vdots & \cdots &\ddots  & \vdots \\
{P} h_{n-1,0} & h_{n-1,1} & \cdots & h_{n-1,k-1} & h_{n-1,k+1} & \cdots & h_{n-1,n-1}
\end{bmatrix}\\
=P\overline{P} \det\begin{bmatrix}
\frac{h_{0,0}}{P\overline{P}} & h_{0,1} & \cdots & h_{0,j-1} & h_{0,j+1} & \cdots& h_{0,n-1} \\
 h_{1,0} & h_{1,1} & \cdots & h_{1,j-1} & h_{1,j+1}& \cdots & h_{1,n-1}  \\
\vdots & \vdots &\ddots & \vdots & \cdots &\ddots  & \vdots \\
h_{k-1,0} & h_{k-1,1} & \cdots & h_{k-1,j-1} & h_{k-1,j+1}& \cdots & h_{k-1,n-1}  \\
 h_{k+1,0} & h_{k+1,1} & \cdots & h_{k+1,j-1} & h_{k+1,j+1}& \cdots & h_{k+1,n-1}  \\
\vdots & \vdots &\ddots & \vdots & \cdots &\ddots  & \vdots \\
 h_{n-1,0} & h_{n-1,1} & \cdots & h_{n-1,j-1} & h_{n-1,j+1} & \cdots & h_{n-1,n-1}
 \end{bmatrix}\\
=r^2 \det\begin{bmatrix}
\frac{h_{0,0}}{r^2} & h_{0,1} & \cdots & h_{0,k-1} & h_{0,k+1} & \cdots& h_{0,n-1} \\
 h_{1,0} & h_{1,1} & \cdots & h_{1,k-1} & h_{1,k+1}& \cdots & h_{1,n-1}  \\
\vdots & \vdots &\ddots & \vdots & \cdots &\ddots  & \vdots \\
h_{j-1,0} & h_{j-1,1} & \cdots & h_{j-1,k-1} & h_{j-1,k+1}& \cdots & h_{j-1,n-1}  \\
 h_{j+1,0} & h_{j+1,1} & \cdots & h_{j+1,k-1} & h_{j+1,k+1}& \cdots & h_{j+1,n-1}  \\
\vdots & \vdots &\ddots & \vdots & \cdots &\ddots  & \vdots \\
 h_{n-1,0} & h_{n-1,1} & \cdots & h_{n-1,k-1} & h_{n-1,k+1} & \cdots & h_{n-1,n-1}
\end{bmatrix}.
\end{align*}
Thus we conclude that $g^{\overline{j},k}$ is irrelevant to $\theta$. We let $h^{j,k}$ be the function irrelevant to $\theta$ satisfying
$$
g^{\overline{j},k} = \begin{cases}
h^{j,k} & \textrm{$(j,k)=(0,0)$ or $jk\neq 0$,} \\
\overline{P}h^{0,k} & \textrm{if $j=0$,} \\
{P}h^{j,0} & \textrm{if $k=0$.} \\
\end{cases}
$$

We let $f=f(Z_1,\cdots,Z_{n-1},1)$. We now show that $\|\nabla W\|^2$ is irrelevant to $\theta$:
\begin{align*}
\|\nabla W\|^2 &= \sum_{j,k} \overline{W_j} g^{\overline{j},k} W_k \\
&=\|f\|^2 g^{\overline{0},0} + \sum_{j,k\neq 0} \overline{W_j} g^{\overline{j},k}W_k + \sum_k \left(P\overline{f}f_k g^{\overline{0},k}+\overline{P}f\overline{f_k}g^{\overline{k},0}\right)\\
&= \|f\|^2 g^{\overline{0},0} +\sum_{j,k\neq 0} \overline{W_j}g^{\overline{j},k}W_k + r^2\sum_k(\overline{f}f_k h^{0,k}+f\overline{f_k}h^{k,0}).
\end{align*}

When $k=0$, $j_0=0$, $\overline{W}_0 g^{\overline{0},0}$ is irrelevant to $\theta$, and
$$
\overline{\partial}(\overline{W_0} g^{\overline{0},0}) =P\frac{\overline{f}}{2r}\frac{\partial g^{\overline{0},0}}{\partial r}d\overline{P}+ \sum_{l=1}^{n-1}\frac{\partial (\overline{f}g^{\overline{0},0})}{\partial\overline{X}_l} d\overline{X}_l.
$$
When $k\neq 0$ and $j_k\neq 0$, $\overline{W_k}g^{\overline{k},j_k} = \overline{P}\overline{f_k} g^{\overline{k},j_k}$, and $\overline{f_k} g^{\overline{k},j_k}$ is irrelevant to $\theta$.
$$
\overline{\partial}(\overline{W_k}g^{\overline{k},j_k}) = \overline{P}\left(\sum_{l=1}^{n-1} \frac{\partial (\overline{f_k} g^{\overline{k},j_k})}{\partial \overline{X}_l}d\overline{X}_l\right)+\overline{f}_k \left(\frac{r}{2}\frac{\partial g^{\overline{k},j_k}}{\partial r} + g^{\overline{k},j_k}\right)d\overline{P}.
$$
When $k=0$ but $j_0\neq 0$, $\overline{W_0}g^{\overline{0},j_0}=\overline{P}\overline{f}h^{0,j_0}$ and
$$
\overline{\partial} (\overline{W_0}g^{\overline{0},j_0})=\overline{P}\left(\sum_{l=1}^{n-1}\frac{\partial(\overline{f}h^{0,j_0})}{\partial\overline{X}_l}d\overline{X}_l\right)+\frac{\overline{f}}{2r}\frac{\partial h^{0,j_0}}{\partial r} d\overline{P}.
$$
When $k\neq 0$ but $j_k=0$, $\overline{W_k}g^{\overline{k},0}=r^2\overline{f_k}h^{k,0}$ and
$$
\overline{\partial} (\overline{W_k}g^{\overline{k},0})=r^2\left(\sum_{l=1}^{n-1}\frac{\partial(\overline{f_k}h^{k,0})}{\partial\overline{X}_l}d\overline{X}_l\right)+P\left(\overline{f_k}h^{k,0} +\overline{f_k}\frac{r}{2}\frac{\partial h^{k,0}}{\partial r}\right)d\overline{P}.
$$
Thus, comparing the degrees of $P$ and $\overline{P}$ in all terms in $A$, we conclude that $A$ can be written as
\begin{align}
A= \sum_{l=1}^{n-1} \overline{P}^{n-2} B_l d\overline{X}_1\wedge \cdots \wedge \widehat{d\overline{X}_{l}}\wedge \cdots \wedge d\overline{X}_{n-1}\wedge d\overline{P} + \overline{P}^{n-1}Cd\overline{X}_1\wedge\cdots\wedge d\overline{X}_{n-1},
\end{align}
where $B_l$ and $C$ are irrelevant to $\theta$. We write 
$$
B:= \sum_{l=1}^{n-1} B_l d\overline{X}_1\wedge \cdots \wedge \widehat{d\overline{X}_{l}}\wedge \cdots \wedge d\overline{X}_{n-1},
$$
and write $C:= Cd\overline{X}_1\wedge\cdots\wedge d\overline{X}_{n-1}$ abusively. Then,
$$
A = \overline{P}^{n-2}B\wedge d\overline{P} + \overline{P}^{n-1} C.
$$

Therefore,
\begin{align*}
& \int_{0\leq \theta\leq 2\pi} P^j (\overline{\partial}\rho) A\wedge \Omega_{CY}\\
&=  \frac{c}{\|\nabla W\|^{2n}}\left(\int_{0\leq \theta\leq 2\pi}\pm P^j \overline{P}^{n-2}d\overline{P}\wedge dP\wedge(\overline{\partial}\rho)\wedge B \pm P^j \overline{P}^{n-1}dP\wedge (\overline{\partial}\rho)\wedge C\right)\wedge dX_1\wedge \ldots \wedge dX_{n-1}.
\end{align*}
In order to get the volume form, the integration inside the parenthesis must contain $dr$ and $d\theta$ terms. Also, In order for the whole value to be nonzero, the value inside the parenthesis must be nonzero:
$$
\int_{0\leq \theta\leq 2\pi}\pm P^j \overline{P}^{n-2}d\overline{P}\wedge dP\wedge (\overline{\partial}\rho)\wedge B \pm P^j \overline{P}^{n-1}dP\wedge (\overline{\partial}\rho)\wedge C.
$$
Note that
\begin{align*}
d\overline{P}\wedge(\overline{\partial}\rho)\wedge B &= d\overline{P} \wedge\left(\sum_{j=1}^{n-1}\frac{\partial \rho}{\partial\overline{X}_j}d\overline{X}_j\right) \wedge B,
\\
(\overline{\partial}\rho)\wedge C &= \frac{1}{2}\frac{\partial \rho}{\partial r}(dr-ird\theta)\wedge C.
\end{align*}
Using the change of variables, we get
\begin{align*}
\int_{0\leq \theta\leq 2\pi} P^j \overline{P}^{n-2}d\overline{P}\wedge dP &= -2ir^{n+j-1}\left(\int_0^{2\pi}e^{i(j-(n-2))\theta} d\theta\right)dr, \\
\int_{0\leq \theta\leq 2\pi}P^j \overline{P}^{n-1}dP\wedge\frac{\partial \rho}{\partial r}(dr-ird\theta) &= \int_{0\leq \theta\leq 2\pi}P^j \overline{P}^{n-1}ire^{i\theta}d\theta\wedge\frac{\partial \rho}{\partial r}dr+\int_{0\leq \theta\leq 2\pi} P^j \overline{P}^{n-1}e^{i\theta}dr\wedge\left(-ir\frac{\partial\rho}{\partial r}d\theta\right)\\
&=2ir^{n+j}\frac{\partial\rho}{\partial r}\left(\int_0^{2\pi}e^{i(j-(n-1)+1)\theta}d\theta \right)\wedge dr.
\end{align*}
Both are nonzero only when $j=n-2$.

(2) We move on to the problem of computation of the integral
$$
\Tr_{LW}(T_{\rho}(\alpha_q))=\int_{X_{CY}}T_{\rho}(\alpha_q\vdash\Omega_{CY})\wedge\Omega_{CY}.
$$
When $q\neq  0$, it is obvious that $\Tr_{LW}(T_{\rho}(\alpha_q))=0$.
Recall that from \eqref{nonprgenerator}, $\alpha_q:= \U_{LW}^{-1}(\omega_0^{\wedge (n/2)}-\alpha\wedge\omega_0^{\wedge(n/2-1)})$ is a generator of the ``non-primitive" part $\bC$ of $H^0(PV(X_{CY}),\overline{\partial}_W)$ when $n$ is even. Let us write $n'=n/2>1$. We need to compute
\begin{align*}
&\int_{X_{CY}} T_\rho\left(\omega_0^{\wedge n'}-\alpha\wedge\omega_0^{\wedge(n'-1)}\right)\wedge \Omega_{CY}
\\&=\int_{X_{CY}}(\overline{\partial}\rho) V_W \left([\overline{\partial},V_W]^{n'-1}(\omega_0^{\wedge n'})-[\overline{\partial},V_W]^{n'}(\alpha\wedge\omega_{0}^{\wedge(n'-1)})\right)\wedge \Omega_{CY}\\
&=\int_{X_{CY}}(\overline{\partial}\rho)V_W[\overline{\partial},V_W]^{n'-1}\left(\omega_0^{\wedge n'}-[\overline{\partial},V_W]\left(\alpha\wedge\omega_0^{\wedge(n'-1)}\right)\right)\wedge\Omega_{CY}.
\end{align*}

\begin{align*}
&\omega_0^{\wedge n'}-[\overline{\partial},V_W](\alpha\wedge\omega_{0}^{\wedge(n'-1)}) \\&=\omega_0^{\wedge n'}-\overline{\partial} V_W(\alpha\wedge\omega_{0}^{\wedge(n'-1)}) -V_W \overline{\partial} (\alpha\wedge\omega_{0}^{\wedge(n'-1)}) \\
&=\omega_0^{\wedge n'}-\overline{\partial} V_W(\alpha\wedge\omega_{0}^{\wedge(n'-1)}) -V_W( dW\wedge \omega_{0}^{\wedge n'})\\
&=\omega_0^{\wedge n'}-\overline{\partial} V_W(\alpha\wedge\omega_{0}^{\wedge(n'-1)}) -  \omega_{0}^{\wedge n'} +dW\wedge V_W(\omega_0^{\wedge n'}) \\
&=-\overline{\partial} V_W(\alpha\wedge\omega_{0}^{\wedge(n'-1)}) +n' dW\wedge V_W(\omega_0)\wedge\omega_0^{\wedge (n'-1)} \\
&=-\overline{\partial} V_W(\alpha)\wedge\omega_{0}^{\wedge(n'-1)}-dW\wedge\omega_0\wedge V_W(\omega_0^{\wedge(n'-1)})-\alpha\wedge (\overline{\partial} V_W(\omega_0^{\wedge(n'-1)})) +n' dW\wedge V_W(\omega_0)\wedge\omega_0^{\wedge (n'-1)}\\
&=-\overline{\partial} V_W(\alpha)\wedge\omega_{0}^{\wedge(n'-1)}-(n'-1)\alpha\wedge (\overline{\partial} V_W(\omega_0))\wedge\omega_0^{\wedge (n'-2)} + dW\wedge V_W(\omega_0)\wedge\omega_0^{\wedge (n'-1)}
\end{align*}
We have 
$$ V_W\overline{\partial}V_W\alpha = -V_W\omega_0+\overline{\partial}V_W\left(\frac{dP}{nP}-\theta_n\right)$$ 
and $\overline{\partial} V_W\overline{\partial}V_W\alpha = -\overline{\partial}V_W\omega_0.$ Inducticely, we have four constants $A_k,B_k,C_k,D_k$ such that
\begin{align*}
&\!\!\!\!\!\![\overline{\partial},V_W]^{k}(\omega_0^{\wedge n'}-[\overline{\partial},V_W](\alpha\wedge\omega_{0}^{\wedge(n'-1)})) \\
=&A_k(\overline{\partial}V_W\alpha)\wedge(\overline{\partial}V_W\omega_0)^{\wedge k}\wedge \omega_0^{\wedge(n'-k-1)}+B_k\alpha\wedge(\overline{\partial}V_W\omega_0)^{\wedge k+1}\wedge \omega_0^{\wedge(n'-k-2)}\\
&+C_kdW\wedge V_W(\omega_0)\wedge(\overline{\partial}V_W\omega_0)^{\wedge k}\wedge \omega_0^{\wedge(n'-k-2)}+D_k (\overline{\partial}V_W\omega_0)^{\wedge k}\wedge \omega_0^{\wedge(n'-k)}.
\end{align*}
When $k=n'-1$, since the result must be $(1,n-1)$-type, $B_{n'-1}=C_{n'-1}=0$, so
\begin{align*}
&[\overline{\partial},V_W]^{n'-1}(\omega_0^{\wedge n'}-[\overline{\partial},V_W](\alpha\wedge\omega_{0}^{\wedge(n'-1)})) \\
&=A_{n'-1}(\overline{\partial}V_W\alpha)\wedge(\overline{\partial}V_W\omega_0)^{\wedge (n'-1)}+D_{n'-1} (\overline{\partial}V_W\omega_0)^{\wedge (n'-1)}\wedge \omega_0,\\
&V_W[\overline{\partial},V_W]^{n'-1}(\omega_0^{\wedge n'}-[\overline{\partial},V_W](\alpha\wedge\omega_{0}^{\wedge(n'-1)})) \\&=(D_{n'-1}-A_{n'-1})V_W(\omega_0)\wedge(\overline{\partial}V_W\omega_0)^{\wedge (n'-1)}+A_{n'-1}\overline{\partial}V_W\left(\frac{dP}{nP}-\theta_n\right)\wedge(\overline{\partial}V_W\omega_0)^{\wedge (n'-1)}.
\end{align*}
We need to show that the following two integrals vanish:
\begin{align}
\int_{U_n}(\overline{\partial}\rho)\wedge V_W(\omega_0)\wedge(\overline{\partial}V_W\omega_0)^{\wedge (n'-1)}\wedge \Omega_{CY}=0,\label{tint}\\
\int_{U_n}(\partial\rho)\wedge\overline{\partial}V_W\left(\frac{dP}{nP}-\theta_n\right)\wedge(\overline{\partial}V_W\omega_0)^{\wedge (n'-1)}\wedge\Omega_{CY}=0.\label{tintt}
\end{align}
This is again obtained by computing the order of $\overline{P}$ for each terms. Let us write 
\begin{align*}
\omega_0 &= \sum_{j,k}w_{j,\overline{k}}dX_j\wedge d\overline{X}_k,\\
\theta_n &= \sum_{j} u_j dX_j.
\end{align*}
Then,
\begin{align*}
V_W\left(\omega_0\right) = \overline{P}\left(\frac{1}{\|\nabla W\|^2}\sum_{\substack{1\leq j\leq n-1\\1\leq k\leq n-1}}\left(\overline{f}h^{0,j}w_{j,\overline{k}}+\sum_{l=1}^{n-1}\overline{f_l}g^{\overline{l},j}w_{j,\overline{k}}\right)d\overline{X}_k\right).
\end{align*}
The differential form inside the parenthesis is irrelevant to $\theta$. Let us define $F_k$ as the function in $\ud{X},\overline{\ud{X}}$, and $r$,
$$
F_k(\ud{X},\overline{\ud{X}},r):=  \frac{1}{\|\nabla W\|^2}\sum_{\substack{1\leq j\leq n-1\\1\leq k\leq n-1}}\left(\overline{f}h^{0,j}+\sum_{l=1}^{n-1}\overline{f_l}g^{\overline{l},j}\right)w_{j,\overline{k}}.
$$
Then,
\begin{align*}
V_W\left(\omega_0\right) &= \overline{P}\left(\sum_{k=1}^{n-1}F_k d\overline{X}_k\right),\\
\overline{\partial}V_W\left(\omega_0\right) &=\overline{P}\left(\sum_{1\leq k<l\leq n-1}\left(\frac{\partial F_l}{d\overline{X}_k}-\frac{\partial F_k}{d\overline{X}_l}\right) d\overline{X}_k\wedge d\overline{X}_l\right)+d\overline{P}\wedge\left(\sum_{k=1}^{n-1}\left(F_k+\frac{r}{2}\frac{\partial F_k}{\partial r}\right) d\overline{X}_k\right).
\end{align*}
Also,
$$
V_W\left(\frac{dP}{nP}-\theta_n\right) = \overline{P}\left(\frac{\overline{f}g^{\overline{0},0}}{nr^2\|\nabla W\|^2}+\sum_{j=1}^{n-1}\frac{\frac{1}{n}\overline{f_j}h^{j,0}-\overline{f}h^{0,j}u_j}{\|\nabla W\|^2}-\sum_{1\leq j,k\leq n-1}\frac{\overline{f_k}g^{\overline{k},j}u_j}{\|\nabla W\|^2}\right).
$$
We let $G$ be the function in the parenthesis
$$
G(\ud{X},\overline{\ud{X}},r):= \frac{\overline{f}g^{\overline{0},0}}{nr^2\|\nabla W\|^2}+\sum_{j=1}^{n-1}\frac{\frac{1}{n}\overline{f_j}h^{j,0}-\overline{f}h^{0,j}u_j}{\|\nabla W\|^2}-\sum_{1\leq j,k\leq n-1}\frac{\overline{f_k}g^{\overline{k},j}u_j}{\|\nabla W\|^2},
$$
which is irrelevant to $\theta$. Then
$$
\overline{\partial}V_W\left(\frac{dP}{nP}-\theta_n\right) = \left(\frac{r}{2}\frac{\partial G}{\partial r}+G\right)d\overline{P} + \overline{P} \sum_{j=1}^{n-1}\frac{\partial G}{\partial\overline{X}_j}d\overline{X}_j.
$$
Thus, we conclude that 
\begin{align*}
(\overline{\partial}\rho)\wedge V_W(\omega_0)\wedge(\overline{\partial}V_W\omega_0)^{\wedge (n'-1)} = \overline{P}^{n'-1}D(\ud{X},\overline{\ud{X}},r)d\overline{P}\wedge dX_1 \wedge \cdots \wedge dX_{n-1},\\
(\partial\rho)\wedge\overline{\partial}V_W\left(\frac{dP}{nP}-\theta_n\right)\wedge(\overline{\partial}V_W\omega_0)^{\wedge (n'-1)}=\overline{P}^{n'-1}E(\ud{X},\overline{\ud{X}},r)d\overline{P}\wedge dX_1 \wedge \cdots \wedge dX_{n-1},
\end{align*}
for some function $E$ and $D$. Thus, both integrals \eqref{tint} and \eqref{tintt} vanishes, when $n'>1$.

(3) The existence of a homogeneous polynomial $h(\ud{x})\in \bC[\ud{x}]_{n(n-2)}$ such that $\Tr_{LW}(T_\rho(p^{n-2}h(\ud{x})))\neq 0$ is guaranteed by the nondegeneracy of $K^{(0)}_{LW}$ (see \cite[Theorem 3.20]{LW}).
\end{proof}

\begin{proposition} \label{proph}
For $-(n-2)\leq q,r\leq n-2$ with $q+r\neq 0$, 
\begin{align}
[\alpha_{q}]\wedge[\alpha_{r}]=0 \label{multiplication3},
\end{align}
and there exists a nonzero constant $c_q$ such that
$$
[\alpha_q] \wedge [\alpha_{-q}] = c_q [p^{n-2}h(\ud{x})],
$$
where $h(\ud{x})$ is given in Lemma \ref{integrationcompute}.
\end{proposition}
\begin{proof}
We first show that for each $q$, there is a nonzero constant $c_q$ and $v_q\in H^0(PV(X_{CY}),\overline{\partial}_W)\cap \Ker \widehat{\Tr}_{LW}$ such that
\begin{align}
\label{multiplication2}
[\alpha_q]\wedge [\alpha_{-q}] =c_q [p^{n-2}h(\ud{x})]+ v_q.
\end{align}
This is directly given by the nondegeneracy of $K_{LW}^{(0)}$ on $H(PV(X_{CY}),\overline{\partial}_W)$. We must have $[\mu]\in H(PV(X_{CY}),\overline{\partial}_W)$ such that
$$
K_{LW}^{(0)}([\alpha_q],[\mu])=\Tr_{LW}(T_\rho(\alpha_q \wedge\mu))\neq 0.
$$
However, by the definition of $K_{LW}^{(0)}$ and the computation in Proposition \ref{integrationcompute}, we must have 
$$
[T_\rho(\alpha_q \wedge\mu)]=[\alpha_q]\wedge [\mu] = c_q' [p^{n-2}h(\ud{x})]+ v
$$
for some $v\in \Ker\widehat{\Tr}_{LW}$ and nonzero constant $c_q'$. Since $[p^{n-2}h(\ud{x})]\in H^0(PV(X_{CY},\overline{\partial}_W)$, we must have
$$
[\alpha_q]\wedge [\mu_{-q}] =c_q' [p^{n-2}h(\ud{x})]+ v_q'
$$
for some $v_q'\in H^0(PV(X_{CY}),\overline{\partial})\cap \Ker\widehat{\Tr}_{LW}$ and $[\mu_{-q}]\in H^{-q}(PV(X_{CY}),\overline{\partial}_W)$. If $q\neq 0$, since the cohomology class $[\alpha_{-q}]$ fully generates $H^{-q}(PV(X_{CY}),\overline{\partial}_W)$, we directly have $[\mu_{-q}]=c''_q[\alpha_{-q}]$ for some nonzero constant $c''_q$. By letting $c'_q/c''_q=c_q$ and $v_q=v'_q/c''_q$, we obtain
$$
[\alpha_q]\wedge [\alpha_{-q}] =c_q [p^{n-2}h(\ud{x})]+ v_q.
$$
When $q=0$, this is given by Proposition \ref{multiplication1}: for all $j>0$ and $g(\ud{x})\in \bC[x]_{jn}$,
$$
[p^j g(\ud{x}) \alpha_0]=0.
$$
Since $\Tr_{LW}(T_\rho(\alpha_0))=0$ by Proposition \ref{integrationcompute} and the nondegeneracy of $K_{LW}^{(0)}$, we must have
$$
K_{LW}^{(0)}([\alpha_0],[\alpha_0])\neq 0.
$$
Thus we conclude that
$$
[\alpha_0]\wedge [\alpha_0] = c_0 [p^{n-2}h(\ud{x})]+ v_0
$$
for $c_0\neq 0$ and $v_0\in H^0(PV(X_{CY}),\overline{\partial}_W)\cap \Ker\widehat{\Tr}_{LW}$.

Now we prove that $[\alpha_q]\wedge[\alpha_r]=0$ if $q+r\neq 0$. Suppose $n$ is odd. Then, $q,r$ are odd, and $q+r$ is even. $[\alpha_q\wedge\alpha_r]$ is an element in $H^{q+r}(PV(X_{CY}),\overline{\partial}_W)$, but from the computation of the cohomology done in Section \ref{sec2}, $H^{q+r}(PV(X_{CY}),\overline{\partial}_W)=0$ unless $q+r=0$. Thus,
$$
[\alpha_q\wedge\alpha_r]=0
$$
unless $q+r=0$.

Suppose $n$ is even and $n> 4$. When $q+r\neq 0$, since $[\alpha_{q+r}]$ generates $H^{q+r}(PV(X_{CY}),\overline{\partial}_W)$, we must have
$$
[\alpha_q\wedge\alpha_r] = C[\alpha_{q+r}]
$$
for some constant $C$. Suppose $C\neq 0$. From \eqref{multiplication2}, we know that $[\alpha_{q+r}\wedge\alpha_{-q-r}]$ is nonzero, so we must have
$$
[\alpha_q\wedge\alpha_r\wedge\alpha_{-q-r}]\neq 0.
$$
From the choice of the representative of $\alpha_q$, we have
$$
\alpha_q = \U^{-1}_{LW}\left(\omega_0^{\wedge(\frac{n+q}{2})}-\alpha\wedge\omega_0^{\wedge(\frac{n+q}{2}-1)}\right).
$$
Thus, the antiholomorphic degree of the Hodge type of $\alpha_q$ is at least $(n+q)/2-1$. Similarly, the antiholomorphic degree of the Hodge type of each term of $\alpha_q\wedge\alpha_r\wedge\alpha_{-q-r}$ is at least
$$
\frac{n+q}{2}-1+\frac{n+r}{2}-1 + \frac{n-q-r}{2}-1 = \frac{3}{2}n-3.
$$
On the other hand, $\alpha_q$ does not contain the factor $d\overline{p}$, and so is $\alpha_q\wedge\alpha_r\wedge\alpha_{-q-r}$. Thus, the degree of the antiholomorphic differential cannot be greater than $n-1$. Hence,
$$
\frac{3}{2}n-3 \leq n-1 \Rightarrow n\leq 4.
$$
Therefore, we get the contradiction.

When $n=4$, since $H^q(PV(X_{CY}),\overline{\partial}_W)=\bC[\alpha_{q}]$ when $q=\pm 2$, there exist two constant $C_1$ and $C_2$ such that
$$
[\alpha_2\wedge\alpha_0] = C_1[\alpha_2], \qquad [\alpha_0\wedge\alpha_{-2}]= C_2[\alpha_{-2}].
$$
If one of $C_1$ or $C_2$ is nonzero,
$$
C_1[\alpha_2\wedge\alpha_{-2}]=[\alpha_2\wedge\alpha_0\wedge\alpha_{-2}]=C_2[\alpha_2\wedge\alpha_{-2}],
$$
so the other is also nonzero. Then,
$$
[\alpha_2\wedge\alpha_0\wedge\alpha_0\wedge\alpha_{-2}] = C_1 C_2[\alpha_2\wedge\alpha_{-2}]\neq 0.
$$
Note that $\alpha_2$ is a polyvector field of Hodge type $(1,3)$ minus that of Hodge type $(0,2)$, $\alpha_0$ is that of Hodge type $(2,2)$ minus that of Hodge type $(1,1)$, and $\alpha_{-2}$ is that of Hodge type $(3,1)$ minus Hodge type $(2,0)$. Thus, $\alpha_2\wedge\alpha_0\wedge\alpha_0\wedge\alpha_{-2}$ must be a polyvector field purely of type $(4,4)$. However, this is impossible since none of $\alpha_2$, $\alpha_0$, and $\alpha_{-2}$ contain the factor $d\overline{p}$. Thus, $C_1=C_2=0$.

Finally, we prove that $v_q$ in \eqref{multiplication2} must be zero. Suppose $v_q\neq 0$. Since $\alpha_0$ and $R(W)_0$ fully generates $H^0(PV(X_{CY}),\overline{\partial}_W)$ as a vector space, we can write
$$
v_q = \sum_{j=0}^{n-3} d_j[p^{j}g_j(x)] + d'[\alpha_0]
$$
for some constants $d_j$ $(j=1,\ldots,n-3)$ and $d'$, and homogeneous polynomials $g_j(x)\in \bC[x]_{nj}$ such that $[p^jg(\ud{x})]\neq 0$. Suppose one of $d_0,\ldots,d_{n-3}$ is nonzero. Let $j\in \{0,\ldots,n-3\}$ be the minimal element so that $d_j$ is nonzero. We have $[p^{n-2-j}g'(x)]\in H^0(PV(X_{CY}),\overline{\partial}_W)$ such that $[p^{n-2-j}g'(x)\cdot p^j g_{j}(x)] = [p^{n-2}h(\ud{x})]$. Then,
$$
[p^{n-j-2}g'(x)\wedge\alpha_q\wedge\alpha_{-q}] = d_j [p^{n-2}h(\ud{x})].
$$
However, by Proposition \ref{multiplication1}, $[p^{n-j-2}g'(x)\wedge\alpha_q\wedge\alpha_{-q}]=0$, so $d_j=0$, which is contradictory.

Suppose $d_0,\cdots,d_{n-3}$ are all zero but $d'\neq 0$. Similarly, we have
$$
[\alpha_0\wedge\alpha_q\wedge\alpha_{-q}] = d'c_0 [p^{n-2}h(\ud{x})] + d'c_0 v_0.
$$
When $q\neq 0$, the left hand side is zero by \eqref{multiplication3}. When $q=0$, we have
$$
[\alpha^{\wedge k}_0] = (d')^{k-1}c_0 [p^{n-2}h(\ud{x})] + (d')^k v_0 \neq 0
$$
for all positive integer $k$. However, this is impossible since for $k$ large enough, we must have $[\alpha^{\wedge k}_0]=0$. Thus, we conclude that $v_q=0$ for each $q=1,\ldots,n-1$.
\end{proof}

\begin{comment}
In summary, the multiplication on $\Omega_{LW}$ is described as follows: for $g_1(\ud{x},p),g_2(\ud{x},p)\in \bC[\ud{x},p]_0$, $a_{k}, b_k\in \bC$ for $k=-n+2,-n+4,\ldots,n-2$,
\begin{align*}
\left([g_1(\ud{x},p)]+\sum_{q}a_{q}[\alpha_q]\right)\wedge\left([g_2(\ud{x},p)]+\sum_{q}b_{q}[\alpha_q]\right)\\
=[g_1(\ud{x},p)\cdot g_2({\ud{x},p})] + \sum_{q}a_{q}b_{-q}c_{q}[p^{n-2}h(\ud{x})]+\sum_q\left(a_q g_2(0,0)+b_qg_1(0,0)[\alpha_q]\right).
\end{align*}
\end{comment}
}

\begin{proof}[Proof of Theorem \ref{theoremtwo}]

We recall the residue map
 \[\mathrm{res}: H^n(\C\BP^{n-1} \setminus V,\C) \to H^{n-2}(V; \C)\] is locally defined by sending $\omega = \frac{df}{f} \wedge \omega_1 + \omega _2$, where $\omega_1,\omega_2$ are holomorphic, to $\omega_1$, so it can be shown to respect the Hodge filtration and induces
 \[\mathrm{res}: F^p H^{n-1}(\C\BP^{n-1} \setminus V,\C) \to F^{p-1}H^{n-2}(V; \C).\]

Let
\begin{align*}
\Omega_{\ud x} = \sum_{i=1}^n (-1)^{i-1} x_i dx_0\wedge \cdots \wedge \widehat {d x_i} \wedge \cdots \wedge dx_n
\end{align*}
Then the Griffiths map $\psi_f: \bigoplus_{k=0}^{n-2} R(f)_{kn} \to H^{n-2}(V;\C)=H^{n-2}(\cA(V),\bar\partial)$ is given by
\begin{align} \label{Gmape}
\psi_f([p^i g(\ud x)]) = \mathrm{res} \bigg(\bigg[\frac{(-1)^i i ! g(\ud x) \Omega_{\ud x}}{f(\ud x)^{i+1}}\bigg]\bigg),\quad i=0, 1, \cdots, n-2,
\end{align}
where $\bigg[\frac{(-1)^i i ! g(\ud x) \Omega_{\ud x}}{f(\ud x)^{i+1}}\bigg]$ represents the cohomology class in $H^{n-1}(\C\BP^{n-1} \setminus V,\C)$. 
By Macaulay's theorem
\be
R(f)_{kn} = 0, \quad k \geq n-1
\ee
so that
\be
\bigoplus_{k=0}^\infty R(f)_{kn} =\bigoplus_{k=0}^{n-2} R(f)_{kn}.
\ee

\begin{proposition}\cite{FTY} \label{ringp}
Let $\U_{BK}$ be the contraction with 
\bea \label{volume}
\Omega_f:=(-1)^{n-2}\psi_f(1).
\eea
Then the map 
\bea \label{pbk}
\Psi_{BK}:=\U_{BK}^{-1} \circ \psi_f: \bigoplus_{k=0}^\infty R(f)_{kn} \to (H^0(\PV(V), \bar\partial), \wedge)
\eea
 is a ring isomorphism.
\end{proposition}

\begin{lemma}\label{hlemmatwo}
Let $c_q h(\ud x)\in \C[\ud x]_{n(n-2)}$ be the polynomial given in Proposition \ref{proph}. Fix $\O_{c_q h} \in \PV^{n-2,n-2}(V)$ such that $ \Psi_{BK} ([c_q h(\ud x)])=[\O_{c_q h}]$.
For each $q = -(n-2), -(n-4), \ldots, n-4, n-2$, we can choose $[\b_q] \in H^q(\PV(V), \bar\partial)$ such that
\be
[\b_q] \wedge [\b_{-q}] = [\O_{c_q h}].
\ee
Moreover, for $-(n-2)\leq q, r \leq n-2$ with $q+r \neq 0$,
\be
[\b_q] \wedge [\b_r] =0
\ee
\end{lemma}
\begin{proof}
 For $q =-(n-2), -(n-4), \ldots, n-4, n-2$, the non-primitive part $H^q_{npr}(\PV(V), \bar\partial)$ is one-dimensional by the Weak Lefschetz theorem and
 the Hodge decomposition $H_{npr}^{2i-n}(\PV(V),\bar\partial) = H_{npr}^{n-1-i,i-1}(\PV(V), \bar\partial)\simeq H^{i-1,i-1}_{npr}(\cA(V), \bar\partial) \simeq \C$ exists.
Thus, there exists 
$\b_q \in \PV^{n-1-i, i-1}(V)$ for $i=1, \cdots, n-1$ (with $q=2i-n$) with the desired properties.
\end{proof}

We have a vector space isomorphism
\begin{eqnarray*}
    \C  \xrightarrow{\sim}  H^{2i-n}_{npr}(\PV(V), \bar\partial), \quad e_i \mapsto [\b_{2i-n}], \quad i=1, \ldots, n-1,
\end{eqnarray*}    
By using this, we can extend \eqref{pbk} to obtain a vector space isomorphism
\be
\Psi_{BK}: \bigoplus_{k=0}^\infty R(f)_{kn}\oplus \C^{n-1} \to H(\PV(V), \bar\partial).
\ee

Define a linear map $\phi$ by
\be
\phi: R(W)_0 \to \bigoplus_{k=0}^\infty R(f)_{kn}, \quad [p^k g(\ud x)] \mapsto [g(\ud x)],
\ee
which can be checked to be a well-defined ring isomorphism;
we extend $\phi$ by using the identity map $\C^{n-1} \to \C^{n-1}$:
\be
\phi: R(W)_0 \oplus \C^{n-1} \to  \bigoplus_{k=0}^\infty R(f)_{kn} \oplus \C^{n-1}
\ee

Finally, we define
\be
\Phi:=\Psi_{BK} \circ \phi \circ\Psi_{LW}^{-1}: H(\PV_{W,\infty}(X_{\mathrm{CY}}), \bar\partial_W) \to H(\PV(V), \bar\partial),
\ee
which fits into the following commutative diagram:
\[\begin{tikzcd}
	{H(\PV_{W,\infty}(X_{\mathrm{CY}}),\bar\partial_W)} & {} & {H(\PV(V),\bar\partial)} \\
	{H(\cA_{W,\infty}(X_{\mathrm{CY}}), \bar\partial+dW)} & {R(W)_0\oplus\C^{n-1}\xrightarrow{\phi} \bigoplus_{k=0}^\infty R(f)_{kn} \oplus \C^{n-1}} & {H(\cA(V),\bar\partial).}
	\arrow["{\U_{LW}}"', from=1-1, to=2-1]
	\arrow["{\Psi_{LW}^{-1}}"', from=1-1, to=2-2]
	\arrow["{\U_{BK}}", from=1-3, to=2-3]
	\arrow["{}", from=2-1, to=2-2]
	\arrow["{}"', from=2-2, to=2-3]
	\arrow["{\Psi_{BK}}"', from=2-2, to=1-3]
\end{tikzcd}\]
Note that $\Phi([\a_q]) = [\b_q]$ for $q=-(n-2), \ldots, n-1$.
Proposition \ref{proph}, Proposition \ref{ringp}, and Lemma \ref{hlemmatwo} imply that
$\Phi$ is a ring isomorphism.

For the proof of \eqref{ip}, note that
for $\a \in \PV_{W,\infty}^{n-i,n-j} (X_{\mathrm{CY}}), \b \in \PV_{W,\infty}^{i,j} (X_{\mathrm{CY}})$
\be
K_{LW}^{(0)}(\a, \b):=\int_{X_{\mathrm{CY}}} \U_{LW}(\a \wedge \b) \wedge \O_{\mathrm{CY}} =(-1)^{nj + (n+1)i} \int_{X_{\mathrm{CY}}} \U_{LW}(\a) \wedge \U_{LW}(\b);
\ee
for $\a \in \PV^{n-2-i,n-2-j} (V), \b \in \PV^{i,j} (V)$
\be
K_{BK}^{(0)}(\a,\b):=\int_{V} \U_{BK}(\a \wedge \b) \wedge \Omega_f =(-1)^{(n-2)j + (n-1)i} \int_{V} \U_{BK}(\a) \wedge \U_{BK}(\b).
\ee
Moreover, we have an induced isomorphism $\Tr_{BK}:H^{n-2,n-2}(\PV(V), \bar\partial)\xrightarrow{\simeq} \C$.

By Proposition \ref{integrationcompute}, we have an isomorphism
$H(\PV_{W,\infty}(X_{\mathrm{CY}}), \bar\partial_W)/\Ker(\Tr_{LW}) \simeq \C$.
Since $\Tr_{BK}, \Tr_{LW}$, and $\Phi$ in the following diagram are isomorphisms, there exists a non-zero complex number $c_{\Phi}$ such that the diagram commutes:
\[\begin{tikzcd}
	{H^{n-2,n-2}(\PV(V), \bar\partial)} & \C \\
	{H(\PV_{W,\infty}(X_{\mathrm{CY}}), \bar\partial_W)/\Ker(\Tr_{LW})} & \C.
	\arrow["{\Tr_{BK}}", from=1-1, to=1-2]
	\arrow["\Phi", from=2-1, to=1-1]
	\arrow["{\Tr_{LW}}", from=2-1, to=2-2]
	\arrow["{c_{\Phi}}"', from=2-2, to=1-2]
\end{tikzcd}\]

Using this, we compute
\be
K_{BK}^{(0)}(\Phi(\a), \Phi(\b))&=& \Tr_{BK}(\Phi(\a) \wedge \Phi(\b)) \\
&=& \Tr_{BK}(\Phi(\a\wedge \b))\\
&=&c_\Phi \Tr_{LW}(\a \wedge \b)\\
&=& c_\Phi K_{LW}^{(0)}(\a, \b),
\ee
which completes the proof of \eqref{ip}.
\end{proof}

\begin{comment}
The map $\Tr_{LW}^{\O_{CY}}:\PV_{W,\infty}(X_{\mathrm{CY}}) \to \C$ induces a map
\begin{eqnarray*}
    \Tr_{LW}: H(\PV_{W,\infty}(X_{\mathrm{CY}}), \bar\partial_W) \to \C
\end{eqnarray*}
     and $\ker(\Tr_{LW})$ has codimension 1, i.e. we have an isomorphism
     \begin{eqnarray*}
          \Tr_{LW}:H(\PV_{W,\infty}(X_{\mathrm{CY}}), \bar\partial_W)/ \ker(\Tr_{LW}) \xrightarrow{\simeq}\C.
     \end{eqnarray*} 

Since $\Tr_{BK}^{\O_f}, \Tr_{LW}^{\O_{\mathrm{CY}}}$, and $\Phi$ in the following diagram are isomorphisms, there exists a non-zero complex number $c_{\Phi}$ such that the diagram commutes:
\[\begin{tikzcd}
	{H^{n-2,n-2}(\PV(V))} & \C \\
	{H^{n,n}(\PV_{W,\infty}(X_{\mathrm{CY}}), \bar\partial_W)  } & \C.
	\arrow["{\Tr_{BK}}", from=1-1, to=1-2]
	\arrow["\Phi", from=2-1, to=1-1]
	\arrow["{\Tr_{LW}}", from=2-1, to=2-2]
	\arrow["{c_{\Phi}}"', from=2-2, to=1-2]
\end{tikzcd}\]

Using this, we compute
\be
K_{BK}^{(0)}(\Phi(\a), \Phi(\b))&=& \Tr_{BK}^{\O_f}(\Phi(\a) \wedge \Phi(\b)) \\
&=& \Tr_{BK}^{\O_f}(\Phi(\a\wedge \b))\\
&=&c_\Phi \Tr_{LW}^{\O_{\mathrm{CY}}}(\a \wedge \b)\\
&=& c_\Phi K_{LW}^{(0)}(\a, \b),
\ee
which completes the proof of \eqref{ip}.
\end{comment}

\subsection{Some remark on Frobenius manifolds}\label{sec3.3}

Let $(\cA, \delta, \Delta)$ be a unital dGBV algebra (\cite[Definition 3.2]{LW}) such that
\begin{enumerate}
\item there is $K$ which is a sesquilinear pairing on $\cA[[u]]$ with respect to which $\delta$ is graded skew-symmetric and $\Delta$ is graded symmetric, and $K$ induces a non-degenerate pairing $K^{(0)}$ on $H(\cA,\delta)$ (see \cite[Lemma 3.17]{LW});
\item the quantum differential graded Lie algebra $(\cA[[u]], \delta+u \Delta, \{ -, -\}_{\Delta})$ is smooth formal (see \cite[page 37]{LW}), where
\be
\{\a,\b\}_\Delta:=\Delta(\a \wedge \b) -\Delta\a \wedge \b - (-1)^{|\a|}\a \wedge \Delta(\b),\quad \a,\b \in \cA[[u]].
\ee
\end{enumerate}

For simplicity of notation, let
\be
\Omega:= H(\cA, \delta), \quad \cH:= H(\cA(( u)), \delta+ u \Delta), \quad \cH^{(0)} := H(\cA[[u]], \delta+ u \Delta),
\ee
where $u$ is a formal parameter. Note that $(\O, [1], K^{(0)})$ is a Frobenius algebra with identity $[1]$.
Let $r^{(0)}:\cH^{(0)} \to \O$ be the natural quotient map (corresponding to $u=0$). 

\begin{definition}
\label{spl}
A map $\sigma:\O \to \cH^{(0)}$ is called a splitting if $\sigma$ is a section to the projection $r^{(0)}:\cH^{(0)} \to \Omega$ and $\sigma$ preserves the pairing
\bea \label{sc}
K(\sigma(\a), \sigma(\b)) = K^{(0)} (\a, \b), \quad \forall \a, \b \in \Omega.
\eea
\end{definition}

The construction of Frobenius manifolds in \cite{LLS} and \cite{BK} can be summarized as follows.
\begin{lemma}\label{uniquefrob}
Let $(\cA, \delta, \Delta)$ be a unital dGBV algebra which satisfies $\mathrm{(1)}$ and $\mathrm{(2)}$ above. 
\begin{itemize}
\item Let $\ud \phi =\{\phi_i =[\varphi_i] \in \O \ | \ \varphi_i \in \cA, \ i=1, \ldots, \mu\}$ be a $\C$-basis of $\O$;
\item let $\G$ be a Maurer-Cartan solution\footnote{we refer to \cite[Section 2]{BK} for the notion of Maurer-Cartan solutions to the differential graded Lie algebra (dgla) and their gauge equivalence classes.} of the quantum dgla $(\cA[[\ud t]]((u)), \delta+u\Delta,\{ -, -\})$, whose existence is guaranteed by smooth-formality of $(\cA[[u]], \delta+u \Delta, \{ -, -\})$;
\item let $\sigma:\O \to \cH^{(0)}$ be a splitting in the sense of Definition \ref{spl}.
\end{itemize}
Then there exists a formal Frobenius manifold\footnote{See \cite[Section 1]{BK}.} $(\O,A_{ij}^{k,\ud \phi, [\G], \sigma}, [1], g_{ij}=K^{(0)}(\phi_i, \phi_j))$, where $A_{ij}^{k,\ud \phi, [\G], \sigma} \in \C[[\ud t]]$ is uniquely determined by $\ud \phi$, the gauge equivalence class $[\G]$ of $\G$,
and $\sigma$, such that the Frobenius algebra $$(\O,A_{ij}^{k,\ud \phi,[\G], \sigma}(\ud 0), [1], g_{ij})$$ obtained by evaluating $\ud t=\ud 0$ is isomorphic to the Frobenius algebra $(\O,[1], K^{(0)})$.
\end{lemma}

 Li-Wen's construction utilizes the dGBV algebra $(\PV_{W,\infty}(X_{\mathrm{CY}}), \bar\partial_W, \partial_{\O_{\mathrm{CY}}})$ and Barannikov-Kontsevich' construction utilizes the dGBV algebra $(\PV(V), \bar\partial, \partial_{\Omega_f})$; both \eqref{LWD} and \eqref{BKD} satisfy conditions $\mathrm{(1)}$ and $\mathrm{(2)}$ above (see \cite[Theorem 3.11, Lemma 3.17]{LW} and \cite[Sections 2 and 4]{BK} for details); there are splittings $\sigma_{LW}:\O_{LW}\to \cH_{LW}^{(0)}$ ([Proposition 3.25]\cite{LW}) and $\sigma_{BK}:\O_{BK} \to \cH_{BK}^{(0)}$ (also apply \cite[$f=0$ case of Proposition 3.25]{LW}). 

By smooth formality of the dgla $(\cA_{LW}[[u]], \bar\partial_{W}+u\partial_{\O_{CY}}, \{-,-\}_{\partial_{\O_{CY}}})$, there exists a versal solution to the associated Maurer-Cartan equation, i.e. 
\begin{eqnarray*}
    \G_{LW}=\sum_{i=1}^\mu\varphi_i^{LW} t_i+ \frac{1}{2}\sum_{i,j} \varphi_{i,j}^{LW}t_i t_j + \cdots \in (\cA_{LW}[[\ud t]][[u]])^0
\end{eqnarray*}   
such that
\begin{eqnarray*}
    \bar\partial_W (\G_{LW}) +\frac{1}{2}\{ \G_{LW}, \G_{LW}\}_{\partial_{\O_{CY}}} =0.
\end{eqnarray*}

This implies that $(\cA_{LW}[[\ud t]][[u]], \bar\partial_{W+\G_{LW}}, u\partial_{\O_{CY}})$ is a dGBV algebra; in particular, 
\begin{eqnarray*}
      (\cA_{LW}[[\ud t]], \bar\partial_{W+\G_{LW}}=\bar\partial+\{W+\G_{LW},- \}_{\partial_{\O_{CY}}})
\end{eqnarray*}
is a CDGA (commutative differential graded algebra); the cohomology $\O_{\G_{LW}}:=H(\cA_{LW}[[\ud t]], \bar\partial_{W+\G_{LW}})$ has an induced algebra structure.

    Let $\ud \phi^{LW}$ be a $\C$-basis of $\O_{LW}=H(\cA_{LW},\bar\partial_W)$.
    We ask whether we can find a Maurer-Cartan solution $\Gamma_{BK}$ associated to the dgla $(\cA_{BK}, \bar\partial, \{-,-\}_{\partial_{\O_f}})$ such that
    \begin{equation*}
        A_{ij}^{k,\ud \phi^{LW}, [\G_{LW}], \sigma_{LW}}=A_{ij}^{k,\Phi(\ud \phi^{LW}), [\G_{BK}], \sigma_{BK}}, \quad  c_\Phi g_{ij}^{LW}=g_{ij}^{BK}.
    \end{equation*} 
    Here $g_{ij}^{LW}=K^{(0)}_{LW}(\phi_i,\phi_j)$ and $g_{ij}^{BK}=K^{(0)}_{BK}(\Phi(\phi_i),\Phi(\phi_j))$.
    Answering this question in an explicit way seems to be difficult, because we cannot find a morphism of dGBV algebras between $(\cA_{LW}, \bar\partial_W, \partial_{\O_{CY}})$ and $(\cA_{BK}, \bar\partial, \partial_{\O_f})$ (we cannot even find a cochain map between $(\cA_{LW}, \bar\partial_W)$ and $(\cA_{BK}, \bar\partial)$ which induces an isomorphism on cohomologies). 

\end{document}